\documentclass[11pt]{article}
 \newcommand{\QED}{\hfill \thicklines \framebox(6.6,6.6)[l]{}}

\newcommand{\pend}{\hfill \thicklines \framebox(5.5,5.5)[l]{}}

\usepackage{fancyhdr,graphicx}
\usepackage{amsfonts}
\usepackage{epstopdf}
\usepackage{amsmath,bm}
\usepackage{epic,curves}
\usepackage{eepic}
\usepackage{eucal}
\usepackage{multirow} 
\usepackage{xcolor}
\usepackage{color}%
\usepackage{cite}%
\usepackage{amsmath}
\usepackage{multirow}%
\usepackage{booktabs}%
\usepackage{float}%
\usepackage{algorithm}
\usepackage{algorithmic}
\usepackage{subfigure}

\newtheorem{theorem}{Theorem}
\newtheorem{lemma}{Lemma}

\begin{document}
\title{\bf Equilibrium and Socially optimal of a double-sided queueing system with two-mass point matching time }
\author{Zhen Wang$^{a}$, Cheryl Yang$^{b}$, Liwei Liu$^{a}$, Yiqiang Q. Zhao$^{b,*}$ \\
{\small\em  $^{a}$School of Science, Nanjing University of Science and Technology, Nanjing 210094, Jiangsu, China}\\
{\small\em  $^{b}$  School of Mathematics and Statistics, Carleton University, 1125 Colonel By Drive, Ottawa, ON K1S 5B6, Canada  }\\
}
\date{}
\footnotetext [1] {E-mail addresses:
yiqiangzhao@cunet.carleton.ca
}
\renewcommand{\thefootnote}{*}
\footnotetext [1] {Corresponding author.}
\maketitle

\begin{abstract}

We study a passenger-taxi double-ended queue with impatient passengers and two-point matching time in this paper. The system considered in this paper is different from those considered in the existing literature, which fully considers the matching time between passengers and taxis, and the taxi capacity of the system. The objective is to get the equilibrium joining strategy and the socially optimal strategy under two information levels. For the practical consideration of the airport terminal scenario, two different information levels are considered. The theoretical results show that the passenger utility function in the partially observable case is monotonic. For the complex form of social welfare function of the partially observable case, we use a split derivation. The equilibrium strategy and socially optimal strategy of the observable case are threshold-type. Furthermore, some representative  numerical scenarios are used to visualize the theoretical results. The numerical scenarios illustrate the influence of parameters on the equilibrium strategy and socially optimal strategy under two information levels. Finally, the optimal social welfare for the two information levels with the same parameters are compared.


\end{abstract}

\section{Introduction}
\label{sec:1}

The transferability of the airport to the arriving passengers is one of the main factors for the performance evaluation of the airport. The taxis waiting capacity at the airport is always limited, and taxis stay in a concentrated area. This area usually has some distance from the terminal, which makes an immediate transportation for passengers not always possible. So there is a random variable matching time between matching-passengers and taxis. As there are many means of transportation at the airport, taking a taxi is not the only choice for arriving passengers. Therefore, it is worth exploring whether the arriving passengers would join the taxi queues with a non-zero matching time. The model considered in this paper is based on observation at the John F. Kennedy International Airport (JFK) in New York City. According to our best knowledge, the airport taxi pick-up scenario has not been studied as a double-ended queue with non-zero matching time. Curry and Vany\cite{Curry1978A} studied an airport pick-up scenario with a single-ended queue and public transport competitions. Passos et al.\cite{2013Passos} compared airport taxi transfer models through simulation. Conway et al.\cite{Conway2018Challenges} summarized the taxi structure of several airports, focusing on the centralized taxi area and scheduling program of the JFK International Airport in New York City. Yazici et al. \cite{Yazici2016Modeling} further elaborated the structure of the JFK Airport, putting forward some policy suggestions by using the method of equilibrium strategy, and analyzed the logistic regression model of taxi decision-making. It is worth noting that these models do not consider the double-sided queues between passengers and taxis. Of course, the airport pick-up scenario can be extended to the ride-hailing apps service scenario. After receiving the order, the driver needs to drive some distance from the current locations to receive the passengers.

Kendall \cite{Kendall1951Some} first introduced the double-sided queuing system under the background of studying passenger-taxi in 1951. He proposed a simple model, i.e., taxis and passengers arrive according to a Poisson process, match, then leave the system instantaneously. Dobbie \cite{Dobbie1961Letter} useed Irwin-Skellam distribution to study the double-sided queuing system, which has the arrival rate varying with time. This model complemented Kendall's research. Wang et al. \cite{wang2020equilibrium} studied a double-sided batch transfer queuing system with a gated policy, in which a batch of passengers is matched with a passenger-car. Giveen \cite{giveen1963taxicab} presented a solution and properties of the time-dependent passenger-taxi model, including the conditions for the existence of the limit distribution. Other double-ended queuing systems can be found in the literature \cite{Conolly2002Double,Kashyap1965A}. These studies on passenger-taxi queuing systems ignored the matching time, i.e., they are transferred immediately after matching.

Considering the strategic behavior of customers in service systems can help the customer's subjective behavior more logically. The joining strategy of customers was first studied by Naor \cite{Naor1969The}, which is the $M/M/1$ system under the observable case. He determined the individual equilibrium strategy and the socially optimal joining strategy. Subsequently, the unobservable case was supplemented by Edelson and Hilderbrand \cite{Edelson1975Congestion}. In recent decades, the research of customer joining strategy and socially optimal joining strategy in queuing systems has been reviewed by many scholars, such as, Jiang et al. who studied service systems with transfers of customers in an alternating environment in \cite{Jiang2019Tail}, catastrophe or clearing systems with equilibrium in \cite{boudali2012optimal,economou2013equilibrium}, and queueing systems with setup or vacation in \cite{burnetas2007equilibrium,sun2017equilibrium}. Bu et al. \cite{Bu2020Strategic} proposed a clearing queueing system with $N$ -policy and stochastic restarting scheme. The monograph of Hassin\cite{Hassin2003To} summarized a large number of related studies on joining strategies and optimal social welfare.

The existing literature dealing with taxi and passenger double-sided queuing models in which leave the system immediately after a successful matching. Shi and Lian \cite{Shi2016Optimization} studied a taxi queuing model of two levels of information, which did not consider the double-sided queue and finite taxi capacity. Wang and Liu \cite{Wang2019Equilibrium} investigated a double-ended queueing system with dynamic taxi arrival rate and finite capacity of taxis, but did not consider the impatient behavior of passengers and the matching time. For double-ended queuing systems with non-zero matching time, the complexity of studies is greatly increased because it can not be simplified as the difference between two random variables. Kim et al. \cite{Kim2010Simulation} studied the double-ended queuing problem with non-zero matching time by simulation. Shi et al. \cite{Shi2015Study} presented a passenger-taxi queueing system with non-zero matching time by using the matrix-analytic method, but their main work was done through computations.

The model considered in this paper is particularly relevant in the scenario of the airport taxi transfer. Taxis wait in a limited capacity waiting area, which is called centralized taxi holding (CTH). The airport taxi dispatcher will direct them to the terminal station to pick up arriving passengers. Taxis from the CTH to the terminal are assigned to passengers according to certain rules. It can be assumed that the dispatcher will dynamically control the taxi to the terminal according to the number of passengers at the terminal. The matching time is random. The following practical scenario can justify the model considered in this paper: a passenger in need of service contacts a driver through the app of a ride-hailing service on his mobile phone, after receiving the order, the driver will pick up the passengers from the current location. Based on the scenario of ride-hailing service, passengers need to wait for a random time to get the service. The model considered in this paper is based on the terminal pick-up scenario, and it fully considers the capacity of the taxis, dynamically controls of the taxis, the matching time between the passenger and the taxi, and the impatient behavior of the passengers. So this model has the following characteristics:
\begin{itemize}
\item The model considered in this paper is based on observation at the John F. Kennedy International Airport (JFK) in New York City and the model assumption can be justified by real scenarios.

\item According to our best knowledge, this model is the first to include the passenger-taxi double-ended queue with impatient passengers and two-point matching time.

\item The two information levels considered are consistent with the degree of information acquisition at the terminal or the actual scenario.

\item The theoretical results show that the passenger utility function in the partially observable case is monotonic. For the complex form of social welfare function of the partially observable case, we use a split derivation.

\item The equilibrium strategy and socially optimal strategy of the observable case are threshold-type.

\item The presentation of numerical experiments on typical cases reveals the necessity of information disclosure or hiding, which helps to achieve optimal social welfare.

\end{itemize}

The remaining paper is organized as follows: In Section \ref{sec:2}, the model description and parameter assumptions are given. The equilibrium strategy and the socially optimal joining strategy for the partially observable case are determined in Section \ref{sec:3}. In Section \ref{sec:4}, it can be seen that equilibrium strategy and socially optimal strategy are given in the form of threshold. In Section \ref{sec:5}, some typical scenarios illustrate the influence of parameters on the equilibrium strategy and socially optimal strategy under two information levels. Section \ref{sec:6} presents discussions and further studies.

\section{Model description}
\label{sec:2}

\begin{figure}[ht]
\centering
\includegraphics[width=0.9\textwidth]{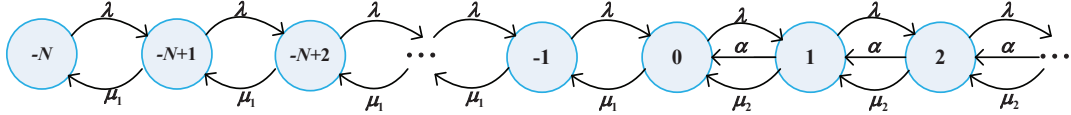}
\caption{The state transition diagram of passenger-taxi double-ended queue with zero matching time.}
\label{Fig:1}
\end{figure}

In this paper, a passenger-taxi double-ended queue system with taxi dynamic controls,  impatience of passengers and two-mass point matching time are considered. The matching time of this system for passengers and taxis is non-zero. In order to consider the strategic behavior of customers, we first determine the performance measure of zero-matching time of this system and introduce the model parameters of zero matching time. The taxis have a capacity of $N$, which means that if there are $N$ taxis in the system, the new taxis will not be able to join the system. So this double-ended queuing system considered in this paper has state space $\Omega = \{ -N, -N+1, -N+2, \cdots, -1, 0, 1, 2, \cdots\}$. $N(t)$ defined as the queue length of taxi or passenger at time $t$. If $N(t) > 0$, the number of passengers in the system is $N(t)$ and there is no taxi queue. If $N(t) < 0$,
the number of taxis in the system is $-N(t)$ and there is no passenger queue. If $N(t) = 0$, there is no taxi nor passenger. Passengers and taxis arrive according to Poisson process.  Passengers (one to four passengers traveling together are considered as one passenger) arrive to the queueing system according to a Poisson process with rate $\lambda$. When there are waiting taxis in the system, passengers directly join the system. When there is no taxi in the system, passengers can choose to join the system with probability $q$, or balk, and the passengers who join the system may have impatient behavior (i.e. they may be exit after joining the system), the impatient time of the passenger is exponentially distributed with rate $\alpha$. Obviously, $\{N(t), t\geq0\}$ is a one-dimensional continuous time Markov chain with state space $\Omega$. The state transition diagram of passenger-taxi double-ended queue with zero matching time is shown in Fig.~\ref{Fig:1}.

The dynamic control of taxi depends on the state $N(t)$ of the system. If there is no passenger (i.e. $N(t) \leq 0$) waiting in the system, the taxi arrival rate is $\mu_{1}$, otherwise (i.e. $N(t) > 0$) the taxi arrival rate is $\mu_{2}$. Obviously, the arrival rate of taxis with passengers is higher than that without passengers, i.e. $\mu_{1} < \mu_{2}$. Passengers and taxis match according to the first-in-first-out discipline, once the matching is successful, the pair is considered to have exited the system. Their matching time occurs outside the system, but it must be considered when considering passenger strategic behavior. When passengers fall into the dilemma of joining or balking, the study of strategic behavior is of great significance. Different types of strategies are used according to the information available to customers when they join the system. This decision frame is based on a natural linear reward-cost structure. We assume that every passenger who completes the traveling need to pay a taxi fee $P$, obtains a reward of $R$. At the same time, every joining passenger incurs a waiting cost $C_{p}$ per unit time associated with service, and a waiting cost $C_{M,P}$ per unit time associated with matching time. Let $C_{T}$ and $C_{M,T}$ represent the taxi's cost per unit time of waiting associated with service and matching time, respectively.

\begin{figure}
\centering
\includegraphics[width=0.7\textwidth]{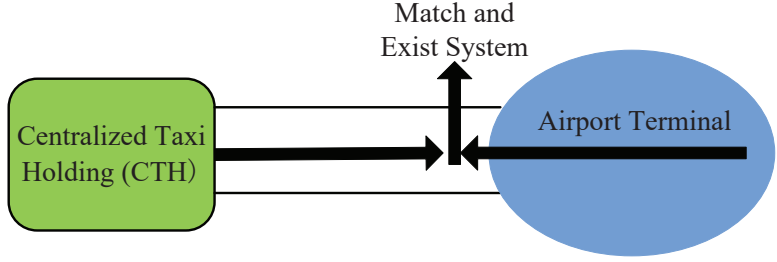}
\caption{Diagram of passenger-taxi double-ended queue in the airport pickup
scenario.}
\label{Fig:2}
\end{figure}

The above-assumed system is particularly relevant in the airport taxi pick-up scenario. In the following sections, we will discuss this system according to the airport taxi pick-up scenario. Taxis wait in a limited capacity waiting area, this area can be called the centralized taxi holding (CTH), the dispatcher is responsible for assigning them to the terminal to receive passengers. It can be assumed that the dispatcher will send taxis to the terminal pickup according to the dynamic control Poisson rates. The matching time, $M $, describes the time it takes for a taxi to get from the CTH to the terminal. Whether the passenger-taxi system with zero-matching time or non-zero matching time, service time refers to the time of loading the passengers and luggage, and send the passengers to the designated place. Passengers arriving at the terminal have many transportation options: taxis, airport shuttle, app-based driving services, and so on. Therefore, the passenger's joining strategy or balk is very important to consider.  Passenger-taxi double-ended queue in the airport pickup scenario is shown in Fig.~\ref{Fig:2}, which shows the taxi queue with finite capacity and the taxi path from the CTH to the terminal, where the taxi matches the passengers in the queue and exits the system.

In the double-ended queueing system of zero matching time described previously, which could be represented as a one-dimensional queue $N(t)$ because the double-ended queueing system of zero matching time implies at least one queue is always empty. The double-ended queueing models with non-zero matching time can be represented as a two-dimensional queueing process $\{N_{P}(t), N_{T}(t)\}$, where $N_{P}(t)$ is the the queues of passengers, $N_{T}(t)$ is the the queues of taxis. Since the double-ended queueing system of non-zero matching time, the queues of passengers and taxis could both be non-zero at once, so it is more complex than the zero matching time system. In this non-zero matching time model, taxis transferred from CTH to the terminal can be considered as artificially assigned to passengers. Therefore, this pair of taxi and passenger can be considered as no longer available in the system, which effectively following the double-ended queueing system with zero matching time. Although taxis and passengers have been artificially ``matched'', there is still a random matching time depending on the state of the system. A realistic scene can be used to describe it, this is similar to the operation of the online car ride-hailing app, where the passenger and driver match, so this pair is not available for the rest of the system, but there is still a waiting time for the driver to pick up the passenger from the current position. Although both passengers and taxis can be considered to be removed from the system within their matching time, the matching time is very important for the strategic behavior of both parties. Therefore, in this model, we can use the simple performance measure of the one-dimensional model while still checking the equilibrium behavior effect of non-zero matching time.

The matching time $M$ between taxis and passengers can be represented by a random variable, which depends on the state of the system, i.e. $N(T)$. We consider a simple two-mass point distribution for this conditional matching time random variable $(M|N(t) = n)$, where $M \in \{ k_1, k_2 \}$ for $k_1 < k_2 \in \mathbb{Z}^+$,
\begin{equation}\label{g01}
P(M = k_1 | N(t) = n) =
\begin{cases}
1, & n \le 0;\\
0 , & \text{otherwise},
\end{cases}
\end{equation}
and
\begin{equation}\label{g02}
P(M = k_2 | N(t) = n) =
\begin{cases}
0, & n \le 0;\\
1 , & \text{otherwise}.
\end{cases}
\end{equation}
It is assumed that this distribution is reasonable because it takes into account the matching time in two cases of the system: when there is no passenger queue, and when there is no taxi queue at the terminal and the passengers are waiting. It stands to reason that when there are no passengers waiting in the terminal, the internal passage between the CTH and the terminal of the airport is unblocked, and the time required for taxis to travel this distance will also be reduced, i.e., $k_{1}< k_{2}$. Based on the above discussion, we have the following two-mass point distribution
\begin{equation}\label{03}
E(M | N(t) = n) =
\begin{cases}
k_1 , & n \le 0;\\
k_2 , & \text{otherwise}.
\end{cases}
\end{equation}
For our analysis to be meaningful, we assume that
\begin{equation}\label{03b}
R>p,
\end{equation}
which ensures that arriving passengers always join the system when they find waiting taxis in the system.

\section{Partially observable case}
\label{sec:3}
\begin{figure} [ht]
\centering
\includegraphics[width=12cm]{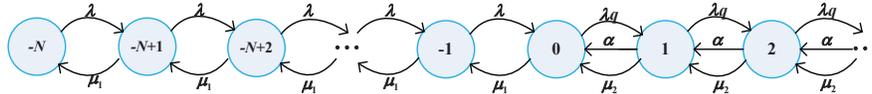}
\caption{The diagram of state transition for partially observable case of the passenger-taxi double-ended queue with zero matching time.}
\label{Fig:g3}
\end{figure}

In this section, we first derive the performance measure of the partially observable case with zero matching time, where the passengers upon arrival can see the taxis queue, but they cannot see the passengers queue (partially observable case). This is similar to the airport terminal, which may show the number of taxis currently in CTH. If there are taxis waiting in the system, then there is no passenger queue, so the arriving passengers directly join the system. Otherwise, the passengers need to wait for new arriving taxis, the arriving passenger joins with probability $q$ and balking with probability $1-q$. Then we derive the passenger equilibrium joining strategy and the socially optimal join strategy of  passenger-taxi double-ended queue system with non-zero matching time according to the performance measure with zero matching time. The state transition diagram with zero matching time is shown in Fig.~\ref{Fig:g3}. For stability, let ${\rho _0}{\rm{ = }}\frac{\lambda }{{{\mu _1}}}< 1$, ${\rho _1}{\rm{ = }}\frac{{\lambda q}}{{\alpha  + {\mu _2}}} < 1$ and ${\rho _2}{\rm{ = }}\frac{\lambda }{{\alpha  + {\mu _2}}}$.

Let ${\pi _n} = \mathop {\lim }\limits_{t \to \infty } P(N(t) = n)$, $(n \in \Omega)$ is the stationary distribution of state $n$ for partially observable case. We derive the stationary distribution by balance equations.
\begin{theorem}  \label{t531}
The stationary distribution of partially observable case of the passenger-taxi double-ended queue
with impatient passengers and zero matching time is given by
\begin{equation}\label{ga1}
{\pi _n} =
\begin{cases}
{\pi _{ - N}}{\rho _0}^{n + N} , &-N \leq n \leq 0;\\
{\pi _{ - N}}{\rho _0}^{N}{\rho _1}^{n}  , & n > 0,
\end{cases}
\end{equation}
where
\begin{equation}\label{ga2}
{\pi _{ - N}} = \frac{{(1 - {\rho _0})(1 - {\rho _1})}}{{1 - {\rho _1} - \rho _0^{N + 1} + \rho _0^N{\rho _1}}}.
\end{equation}
\end{theorem}
{\bf Proof}~~
According to Fig. \ref{Fig:g3}, the following equilibrium equations are obtained
\begin{align}\label{g1}
\lambda {\pi _{ - N}} &= {\mu _1}{\pi _{ - N + 1}},\\
\label{g2}
(\lambda  + {\mu _1}){\pi _n} &= \lambda {\pi _{n - 1}} + {\mu _1}{\pi _{n + 1}}, ~~- N + 1 \le n \le 1,\\
\label{g3}
(\lambda q + {\mu _1}){\pi _0} &= \lambda {\pi _{ - 1}} + (\alpha  + {\mu _2}){\pi _1},\\
\label{g4}
(\lambda q + \alpha  + {\mu _2}){\pi _n} &= \lambda q{\pi _{n - 1}} + (\alpha  + {\mu _2}){\pi _{n + 1}},~~n \ge 1.
\end{align}
From (\ref{g1})-(\ref{g4}) and after some algebraic calculations, we can get (\ref{ga1}). $\pi _{ - N}$ can be obtained by the normalization condition $\sum\limits_{n =  - N}^{ + \infty } {{\pi _n} = 1} $. \QED

By (\ref{ga1}) and (\ref{ga2}), for partially observable case, the effective arrival rate of passengers is
\begin{align}\label{g5}
    \nonumber \lambda _P^* &= \lambda \sum\limits_{n =  - N}^{ - 1} {{\pi _n} + \lambda q\sum\limits_{n = 0}^\infty  {{\pi _n}} }\\
    \nonumber  &= \lambda {\pi _{ - N}}\rho _0^N \left( \sum\limits_{n =  - N}^{ - 1} {\rho _0^n}  + q\sum\limits_{n = 0}^\infty  {\rho _1^n} \right) \\
 &= \lambda {\pi _{ - N}} \left(\frac{{1 - \rho _0^N}}{{1 - {\rho _0}}} + \frac{{q\rho _0^N}}{{1 - {\rho _1}}}\right),
\end{align}
and the effective arrival rate of taxis is
\begin{align}\label{g6}
    \nonumber \lambda _T^* &= {\mu _1}\sum\limits_{n =  - N + 1}^0 {{\pi _n} + {\mu _2}\sum\limits_{n = 1}^\infty  {{\pi _n}} }\\
    \nonumber  & = {\pi _{ - N}}\rho _0^N \left({\mu _1}\sum\limits_{n =  - N + 1}^0 {\rho _0^n}  + {\mu _2}\sum\limits_{n = 1}^\infty  {\rho _1^n} \right)\\
&= {\pi _{ - N}}\left(\frac{{\lambda (1 - \rho _0^N)}}{{1 - {\rho _0}}} + \frac{{{\mu _2}{\rho _1}\rho _0^N}}{{1 - {\rho _1}}}\right).
\end{align}
Denoting $E(L_{P})$ and $E(L_{T})$ as expected queue lengths for passenger and taxi, respectively, they have the following forms
\begin{equation}\label{g7}
E({L_P}) = \sum\limits_{n = 0}^\infty  {n{\pi _n} = {\pi _{ - N}}\rho _0^N\sum\limits_{n = 0}^\infty  {n\rho _1^n} }  = \frac{{{\pi _{ - N}}\rho _0^N{\rho _1}}}{{{{(1 - {\rho _1})}^2}}},
\end{equation}
\begin{equation}\label{g8}
E({L_T}) = \sum\limits_{n =  - N}^0 { - n{\pi _n} = {\pi _{ - N}}\rho _0^N\sum\limits_{n = 0}^N {n{{\left(\frac{1}{{{\rho _0}}}\right)}^n}} }  = \frac{{{\pi _{ - N}}(\rho _0^{N + 1} - {\rho _0} + N - N{\rho _0})}}{{{{(1 - {\rho _0})}^2}}}.
\end{equation}
Using (\ref{g5})-(\ref{g8}) and Little's law, the expected waiting times for passenger and taxi queues can be obtained as follows, respectively,
\begin{equation}\label{g9}
E({W_P}) = \frac{{E({L_P})}}{{\lambda _P^*}} = \frac{{{\pi _{ - N}}\rho _0^N{\rho _1}}}{{\lambda _P^*{{(1 - {\rho _1})}^2}}},
\end{equation}
\begin{equation}\label{g10}
E({W_T}) = \frac{{E({L_T})}}{{\lambda _T^*}} = \frac{{{\pi _{ - N}}(\rho _0^{N + 1} - {\rho _0} + N - N{\rho _0})}}{{\lambda _T^*{{(1 - {\rho _0})}^2}}}.
\end{equation}
\subsection{Equilibrium joining strategy in the partially observable case with two-mass point matching time}
For the partially observable case, we use performance measures of the passenger-taxi double-ended queue with impatient passengers and zero matching time to explore the equilibrium join strategy for the passenger-taxi double-ended queueing with impatient passengers and two-mass point matching time. As mentioned above, after the passenger and taxi are matched successfully, this pair can be regarded as unavailable to the original system, and their matching time is only between this pair. During this matching time, both passenger and taxi of this pair can be considered to be removed from the system. But the matching time between this pair is very important for the study of their equilibrium behavior. We first give the expected waiting time of a joining passenger in this system with zero matching time and no waiting taxi queue (i.e. $N(t) \geq 0$).
\begin{theorem}  \label{t532}
For partially observable case, when the queue length of taxis is zero, the expected waiting time of a joining passenger in the passenger-taxi double-ended queue with impatient passengers and  zero matching time is given by
\begin{equation}\label{g11}
E(W) = \frac{1}{{(\alpha  + {\mu _2}) - \lambda q}}.
\end{equation}
\end{theorem}
{\bf Proof}~~
By (\ref{g7}) and Little's law,
the expected waiting time of a joining passenger in the passenger-taxi double-ended queue with impatient passengers and  zero matching time is
\begin{align}\label{g12}
    \nonumber E(W) &= E(\left. {{W_P}} \right|N(t) \ge 0)\\
    \nonumber  &= \frac{{E({W_P},N(t) \ge 0)}}{{P(N(t) \ge 0)}}\\
    \nonumber &= \frac{{E({L_P})}}{{\lambda q}}/P(N(t) \ge 0)\\
&= \frac{1}{{(\alpha  + {\mu _2}) - \lambda q}},
\end{align}
i.e. (\ref{g11}). \QED

In the partially observable case, the passengers upon arrival can see the taxis queue, but they cannot see the passengers queue. If there are idle taxis in the system, the arriving passengers will not hesitate to join the system. Otherwise, there are no waiting taxis and some passengers queuing in the system, the arriving passenger joins the system with probability $q$. We give the passenger utility function for the partially observable of the passenger-taxi double-ended queue with impatient passengers and two-point matching time by the following lemma \ref{lem1}.
\begin{lemma}  \label{lem1}
For partially observable case, the utility function of an arriving passenger joining the passenger-taxi double-ended queue with impatient passenger and two-point matching time is given by
\begin{equation}\label{g14}
U(q) = R - P - \frac{{{C_P}}}{{(\alpha  + {\mu _2}) - \lambda q}} - {C_{M,P}}E(M),
\end{equation}
where $E(M)$ is the expected matching time for two-mass point distribution,
\begin{equation}\label{g15}
E(M) = {\pi _{ - N}}\rho _0^N \left [{k_1}\frac{{\lambda (1 - \rho _0^{ - N - 1})}}{{\lambda  - {\mu _1}}} + {k_2}\frac{{\lambda q}}{{(\alpha  + {\mu _2}) - \lambda q}} \right ].
\end{equation}
\end{lemma}
{\bf Proof}~~
Based on the reward-cost structure, the utility function of a arriving passenger joining the passenger-taxi double-ended queue with impatient passenger and zero matching time is
\begin{align}\label{g16}
\nonumber {U_1}(q) &= R - P - {C_P}E(W) \\
&= R - P - \frac{{{C_P}}}{{(\alpha  + {\mu _2}) - \lambda q}}.
\end{align}
Because a pair of passenger and taxi match successfully, it is not available to the original system, which corresponding to this pair exiting the system,  and their matching time only occurs between them. Since the matching time is a random variable depending on the state of the system, which is unknown to the arriving passengers, so it does not affect the value of $E(W)$. Let $C_{M,P}$ be the waiting cost per unit time associated with matching time, and $E(M)$ be the expected matching time. So the expected cost incurred due to the matching time is $C_{M,P}E(M)$, where
\begin{align}\label{g17}
    \nonumber E(M) &= \sum\limits_{n =  - N}^\infty  {E(M\left| {N(t) = n} \right.)} {\pi _n}\\
    \nonumber  &= {\pi _{ - N}}\rho _0^N \Big ({k_1}\sum\limits_{n =  - N}^0 {\rho _0^n + {k_2}\sum\limits_{n = 1}^\infty  {\rho _1^n} } \Big )\\
&= {\pi _{ - N}}\rho _0^N\left [{k_1}\frac{{\lambda (1 - \rho _0^{ - N - 1})}}{{\lambda  - {\mu _1}}} + {k_2}\frac{{\lambda q}}{{(\alpha  + {\mu _2}) - \lambda q}}\right ],
\end{align}
i.e. (\ref{g15}). Hence, the utility function of a arriving passenger joining the passenger-taxi double-ended queue with impatient passenger and two-point matching time is
\begin{equation}\label{g18}
U(q) = R - P - \frac{{{C_P}}}{{(\alpha  + {\mu _2}) - \lambda q}} - {C_{M,P}}E(M),
\end{equation}
i.e. (\ref{g14}). \QED

If an arriving passenger finds that there is no taxi queuing in the system, and the arriving passenger needs to face the dilemma of joining or balking. The equilibrium strategy for an arriving passenger is represented by $q_{e}$, $q_{e}$ is given by the following theorem \ref{t533}.
\begin{theorem} \label{t533}
For partially observable case, the equilibrium joining probability $q_{e}$ of a arriving passenger, upon joining the passenger-taxi double-ended queue with impatient passenger and two-point matching time is given
\begin{equation}\label{g19}
    q_{e}=
\begin{cases}
0, & {\rm if}~~ 0< R-P< L_{PO};\\
q_{1}, &{\rm if}~~ L_{PO} \leq R-P \leq V_{PO};\\
1, &{\rm if}~~R-P > V_{PO},
\end{cases}
\end{equation}
where $q_{1}$ is the unique solution of $U(q)=0$ for $q \in [0,1]$, and
\begin{equation}\label{g20}
{L_{PO}} = \frac{{{C_P}}}{{(\alpha  + {\mu _2})}} + \frac{{{C_{M,P}}{k_1}{\mu _1}(1 - {\rho _0})}}{{{\mu _1} - \lambda }},
\end{equation}
\begin{align}\label{g21}
    \nonumber {V_{PO}} &= \frac{{{C_P}}}{{(\alpha  + {\mu _2}) - \lambda }} + \frac{{{C_{M,P}}{k_1}{\mu _1}(1 - \rho _0^{N + 1})(1 - {\rho _0})(\alpha  + {\mu _2} - \lambda )}}{{({\mu _1} - \lambda )[(\alpha  + {\mu _2})(1 - \rho _0^{N + 1}) - \lambda (1 + \rho _0^N)]}} \\
 &+ \frac{{{C_{M,P}}{k_2}\lambda \rho _0^N(1 - {\rho _0})}}{{(\alpha  + {\mu _2})(1 - \rho _0^{N + 1}) - \lambda (1 + \rho _0^N)}}.
\end{align}
\end{theorem}
{\bf Proof}~~
By lemma \ref{lem1}, we can obtain
\begin{align}\label{g22}
    \nonumber U'(q) &=  - \frac{{{C_P}\lambda }}{{{{(\alpha  + {\mu _2} - \lambda q)}^2}}} - \frac{{{C_{M,P}}{k_1}{\lambda ^2}(\alpha  + {\mu _2})\rho _0^{N - 1}(1 - \rho _0^{N + 1})(1 - {\rho _0})(1 + {\rho _0})}}{{({\mu _1} - \lambda ){{[(\alpha  + {\mu _2})(1 - \rho _0^{N + 1}) - \lambda q(1 + \rho _0^N)]}^2}}}\\
& - \frac{{{C_{M,P}}{k_2}\lambda \rho _0^N(1 - {\rho _0})(\alpha  + {\mu _2})(1 - \rho _0^{N + 1})}}{{{{[(\alpha  + {\mu _2})(1 - \rho _0^{N + 1}) - \lambda q(1 + \rho _0^N)]}^2}}},
\end{align}
So $U'(q)<0$ for $q \in [0,1]$. Hence, $U(q)$ is a strictly decreasing function for $q \in [0,1]$, it has unique maximum
\begin{equation}\label{g23}
U(0)=R-P-L_{PO},
\end{equation}
where
\begin{equation*}
{L_{PO}} = \frac{{{C_P}}}{{(\alpha  + {\mu _2})}} + \frac{{{C_{M,P}}{k_1}{\mu _1}(1 - {\rho _0})}}{{{\mu _1} - \lambda }},
\end{equation*}
and it has a unique minimum
\begin{equation}\label{g24}
U(1)=R-P-V_{PO},
\end{equation}
where
\begin{align}
    \nonumber {V_{PO}} &= \frac{{{C_P}}}{{(\alpha  + {\mu _2}) - \lambda }} + \frac{{{C_{M,P}}{k_1}{\mu _1}(1 - \rho _0^{N + 1})(1 - {\rho _0})(\alpha  + {\mu _2} - \lambda )}}{{({\mu _1} - \lambda )[(\alpha  + {\mu _2})(1 - \rho _0^{N + 1}) - \lambda (1 + \rho _0^N)]}} \\
\nonumber  &+ \frac{{{C_{M,P}}{k_2}\lambda \rho _0^N(1 - {\rho _0})}}{{(\alpha  + {\mu _2})(1 - \rho _0^{N + 1}) - \lambda (1 + \rho _0^N)}}.
\end{align}
Therefore, when $R-P \in (0, L_{PO})$, $U(q)$ is negative for every $q \in [0, 1]$, so the best response is balking of $q_{e}=0$.

When $R-P \in [L_{PO}, V_{PO}]$, which indicates that $U(0)>0$ and $U(1)<0$, since the monotonicity of $U(q)$ for $q \in [0, 1]$, there exists a unique solution of the equation $U(q)=0$, let $q_{1}$ be the unique solution. Hence, $q_{1}$ is a proper mixed strategy in this case.

When $R-P \in [V_{PO}, \infty ]$, i.e. $R-P-V_{PO} >0$, it implies $U(q)$ is positive for every $q$. A arriving passenger's best response is always joining, i.e. $q_{e}=1$. \QED
\subsection{Socially Optimal Strategies of Passengers}
Social welfare (social net benefit) is the sum of the welfare of passenger and taxi, which is defined by $S(q)$,
\begin{equation}\label{g25}
S(q) = \lambda _p^*(R - P - {C_P}E({W_P}) - {C_{M,P}}E(M)) + \lambda _T^*(P - {C_T}E({W_T}) - {C_{M,T}}E(M)),
\end{equation}
where $S(q)$ can be written as
\begin{equation}\label{g26}
S(q) = \lambda _p^*(R - P - {C_P}E({W_P})) + \lambda _T^*(P - {C_T}E({W_T})) + E(M)( - \lambda _p^*{C_{M,P}} - \lambda _T^*{C_{M,T}}),
\end{equation}
let $S_{0}(q)=\lambda _p^*(R - P - {C_P}E({W_P})) + \lambda _T^*(P - {C_T}E({W_T}))$ and $S_{M}(q)=E(M)( - \lambda _p^*{C_{M,P}} - \lambda _T^*{C_{M,T}})$, i.e., $S(q)=S_{0}(q)+S_{M}(q)$. Obviously, $S_{0}(q)$ represents the social welfare of the system without matching time, $S_{M}(q)$ is social welfare associated with matching time. From (\ref{g5}) and (\ref{g6}), we can get
\begin{equation}\label{g27}
\lambda _T^* = \lambda _P^* - \frac{{\alpha \rho _0^N{\rho _1}}}{{1 - {\rho _1}}}{\pi _{ - N}},
\end{equation}
so
\begin{align}\label{g28}
    \nonumber {S_0}(q) &= \lambda _p^*(R - P - {C_P}E({W_P})) + (\lambda _P^* - \frac{{\alpha \rho _0^N{\rho _1}}}{{1 - {\rho _1}}}{\pi _{ - N}})(P - {C_T}E({W_T}))\\
&= \lambda _p^*(R - P - {C_P}E({W_P})) + \lambda _P^*(P - {C_T}E({W_T})) - \frac{{\alpha \rho _0^N{\rho _1}}}{{1 - {\rho _1}}}{\pi _{ - N}}(P - {C_T}E({W_T})).
\end{align}
Let
\begin{equation}\label{g29}
{S_1}(q) = \lambda _p^*(R - P - {C_P}E({W_P})) + \lambda _P^*(P - {C_T}E({W_T})),
\end{equation}
and
\begin{equation}\label{g30}
{S_2}(q) =  - \frac{{\alpha \rho _0^N{\rho _1}}}{{1 - {\rho _1}}}{\pi _{ - N}}(P - {C_T}E({W_T}))
\end{equation}
so
\begin{equation}\label{g31}
S(q) = S_{0}(q)+S_{M}(q)={S_1}(q)+{S_2}(q)+S_{M}(q).
\end{equation}

The study of socially optimal joining strategy can promote the maximization of social welfare, the socially optimal joining strategy is represented by $q^{*}$ in the partially observable case. We give the socially optimal join strategy by the following theorem \ref{t534}.
\begin{theorem} \label{t534}
Consider the partially observable case, $S(q)$ is a strictly decreasing function for $q \in [0, 1]$. Then, the socially optimal strategy is ${q^*} = \min (\frac{{dS(q)}}{{dq}}\left| {_{q = {q^*}}} \right. = 0, 1)$.
\end{theorem}
{\bf Proof}~~
From (\ref{g31}), which can get $S(q) ={S_1}(q)+{S_2}(q)+S_{M}(q)$. It is intractable for us to directly calculate the first derivative of $S(q)$ with respect to $q \in [0, 1]$, but we can take the function $S(q)$ apart, then calculate its first derivative of $q$. Obviously,
\begin{equation}\label{g32}
S'(q) = {S_1}^\prime (q) + {S_2}^\prime (q) + {S_M}^\prime (q).
\end{equation}

We first calculate $\frac{{d{S_1}(q)}}{{dq}}$. From (\ref{g29}), (\ref{g9}) and (\ref{g10}), we have
\begin{equation}\label{g33}
{S_1}(q) = \lambda _p^*R - {C_P}E({L_P}) - {C_T}E({L_T}),
\end{equation}
by using (\ref{g5}), (\ref{g7}) and (\ref{g8}), the detailed expression of ${S_1}(q)$ is as follows
\begin{align}\label{g34}
\nonumber {S_1}(q) &= R\frac{{\lambda (1 - \rho _0^N)(1 - {\rho _1}) + \lambda q\rho _0^N(1 - {\rho _0})}}{{1 - {\rho _1} - \rho _0^{N + 1} + \rho _0^N{\rho _1}}} - {C_P}\frac{{\rho _0^N{\rho _1}(1 - {\rho _0})}}{{(1 - {\rho _1})(1 - {\rho _1} - \rho _0^{N + 1} + \rho _0^N{\rho _1})}}\\
&- {C_T}\frac{{(1 - {\rho _1})}}{{(1 - {\rho _0})(1 - {\rho _1} - \rho _0^{N + 1} + \rho _0^N{\rho _1})}},
\end{align}
let $\overline C  = {C_T}(\rho _0^{N + 1} - {\rho _0} + N - N{\rho _0})/{(1 - {\rho _0})^2}$, ${S_1}(q)$ can be simplified as
\begin{equation}\label{g35}
{S_1}(q) = {\mu _2}R - {\pi _{ - N}}\left [ \left (\frac{{{\mu _2}}}{{1 - {\rho _0}}} - \frac{{{\mu _2}\rho _0^{N + 1}}}{{1 - {\rho _0}}} - \frac{{\lambda (1 - \rho _0^N)}}{{1 - {\rho _0}}} \right )R + {C_P}\frac{{\rho _0^N{\rho _1}}}{{{{(1 - {\rho _1})}^2}}} + \overline C \right ].
\end{equation}
The first derivative of ${S_1}(q)$ of $q \in [0,1]$ is
\begin{align}\label{g36}
\nonumber {S_1}^\prime (q) &= \frac{{\rho _0^N{\rho _2}}}{{{{(1 - {\rho _1})}^2}}}\pi _{ - N}^2 \left [ R\frac{{{\mu _2}(1 - \rho _0^{N + 1}) - \lambda (1 - \rho _0^N)}}{{1 - {\rho _0}}} + {C_P}\frac{{\rho _0^N{\rho _1}}}{{{{(1 - {\rho _1})}^2}}} + \overline C  - \frac{{1 + {\rho _1}}}{{1 - {\rho _1}}}\frac{{{C_P}}}{{{\pi _{ - N}}}} \right ]\\
\nonumber &= \frac{{\rho _0^N{\rho _2}\pi _{ - N}^2}}{{{{(1 - {\rho _1})}^4}}} \left [R \frac{{{\mu _2}(1 - \rho _0^{N + 1}) - \lambda (1 - \rho _0^N)}}{{1 - {\rho _0}}}{(1 - {\rho _1})^2} + {C_P}\rho _0^N{\rho _1} + \overline C {(1 - {\rho _1})^2} - \frac{{{C_P}(1 + {\rho _1})(1 - {\rho _1})}}{{{\pi _{ - N}}}} \right ]\\
\nonumber &= \frac{{\rho _0^N{\rho _2}\pi _{ - N}^2}}{{{{(1 - {\rho _1})}^4}}} \left [\left (R\frac{{{\mu _2}(1 - \rho _0^{N + 1}) - \lambda (1 - \rho _0^N)}}{{1 - {\rho _0}}} + \overline C  + {C_P}\frac{{1 - \rho _0^N}}{{1 - {\rho _0}}} \right )\rho _1^2 \right.\\
&- \left. 2{\rho _1}\left (R\frac{{{\mu _2}(1 - \rho _0^{N + 1}) - \lambda (1 - \rho _0^N)}}{{1 - {\rho _0}}} + \overline C \right ) + \left (R\frac{{{\mu _2}(1 - \rho _0^{N + 1}) - \lambda (1 - \rho _0^N)}}{{1 - {\rho _0}}} + \overline C  - {C_P}\frac{{1 - \rho _0^{N + 1}}}{{1 - {\rho _0}}}\right )  \right ].
\end{align}
Since
\begin{equation}\label{g37}
\Delta  = 4{C_P}\rho _0^N \left(\frac{{R{\mu _2}(1 - \rho _0^{N + 1}) - \lambda R(1 - \rho _0^N)}}{{1 - {\rho _0}}} + \overline C \right) + 4C_P^2\frac{{(1 - \rho _0^N)(1 - \rho _0^{N + 1})}}{{{{(1 - {\rho _0})}^2}}} > 4C_P^2\frac{{{{(1 - \rho _0^N)}^2}}}{{{{(1 - {\rho _0})}^2}}},
\end{equation}
we obtain that the larger root $\overline{q}$ of the equation ${S_1}^\prime (q)=0$ is
\begin{equation}\label{g38}
\overline q  = \frac{{ - B + \sqrt \Delta  }}{{2A}} > \frac{{ - B + 2{C_P}(1 - \rho _0^N)/(1 - {\rho _0})}}{{2A}} = 1,
\end{equation}
where
\begin{equation*}
A = \frac{{R{\mu _2}(1 - \rho _0^{N + 1}) - \lambda R(1 - \rho _0^N)}}{{1 - {\rho _0}}} + \overline C  + {C_P}\frac{{1 - \rho _0^N}}{{1 - {\rho _0}}},
\end{equation*}
and
\begin{equation*}
B =  - 2\left(\frac{{R{\mu _2}(1 - \rho _0^{N + 1}) - \lambda R(1 - \rho _0^N)}}{{1 - {\rho _0}}} + \overline C \right).
\end{equation*}
Therefore, since $\overline{q}>1$, which leads to $\frac{{d{S_1}(q)}}{{dq}}<0$ for $q \in [0, 1]$.

Secondly, calculating $\frac{{d{S_2}(q)}}{{dq}}$. By using (\ref{ga2}), (\ref{g10}) and (\ref{g30}), the detailed expression of ${S_2}(q)$ is as follows
\begin{align}\label{g39}
\nonumber {S_2}(q) &=  - \frac{{\alpha \rho _0^N{\rho _1}}}{{1 - {\rho _1}}}{\pi _{ - N}}\left(P - {C_T}E({W_T})\right)\\
&=  - \frac{{\alpha \rho _0^N{\rho _1}(1 - {\rho _0})P}}{{1 - {\rho _1} - \rho _0^{N + 1} + \rho _0^N{\rho _1}}} - \frac{{{C_T}\alpha \rho _0^N(\rho _0^{N + 1} - {\rho _0} + N - N{\rho _0})(1 - {\rho _1}){\rho _1}}}{{(1 - {\rho _1} - \rho _0^{N + 1} + \rho _0^N{\rho _1})[\lambda (1 - \rho _0^N)(1 - {\rho _1}) + {\mu _2}\rho _0^N(1 - {\rho _0}){\rho _1}]}},
\end{align}
The first derivative of ${S_2}(q)$ of $q \in [0,1]$ is
\begin{equation}\label{g40}
{S_2}^\prime (q) =  - P\frac{{{D_1}}}{{{{(1 - {\rho _1} - \rho _0^{N + 1} + \rho _0^N{\rho _1})}^2}}} - {C_T}\frac{{{D_2}{D_3}}}{{{{(1 - {\rho _1} - \rho _0^{N + 1} + \rho _0^N{\rho _1})}^2}{D_4}}},
\end{equation}
where
\begin{align}\label{g41}
\nonumber {D_1} &= \alpha \rho _0^N(1 - {\rho _0}){\rho _2}[(1 - {\rho _1} - \rho _0^{N + 1} + \rho _0^N{\rho _1}) + (1 - {\rho _0})(1 - \rho _0^N)],\\
\nonumber {D_2} &= \alpha \rho _0^N(\rho _0^{N + 1} - {\rho _0} + N - N{\rho _0}),\\
\nonumber {D_3} &= \lambda (1 - \rho _0^N){\rho _2}{(1 - {\rho _1})^2}(1 - \rho _0^{N + 1}) + {\mu _2}\rho _0^{2N}{(1 - {\rho _0})^2}{\rho _2}\rho _1^2,\\
\nonumber {D_4} &= \lambda (1 - \rho _0^N)(1 - {\rho _1}) + {\mu _2}\rho _0^N(1 - {\rho _0}){\rho _1}.
\end{align}
According to the previous hypothesis $\rho_{0}<1$, i.e. $\lambda < \mu_{1}$, we can obtain $D_{1}$, $D_{2}$, $D_{3}$, $D_{4}>0$, so $\frac{{d{S_2}(q)}}{{dq}}<0$ for $q\in [0, 1]$.

Thirdly, calculating $\frac{{d{S_M}(q)}}{{dq}}$. By using (\ref{g5}), (\ref{g6}) and (\ref{g17}), the detailed expression of ${S'_M}(q)$ is as follows
\begin{equation}\label{g42}
{S'_M}(q) =  - \frac{{[{A_1} + {A_3}{B_1}({\mu _1} - \lambda )]({C_{M,T}} + {C_{M,P}}) + [{A_2}{B_1} + {A_4}({\mu _1} - \lambda )][(\alpha  + {\mu _2}){C_{M,P}} + {\mu _2}{C_{M,T}}]}}{{({\mu _1} - \lambda )(1 - {\rho _1} - \rho _0^{N + 1} + \rho _0^N{\rho _1})}},
\end{equation}
where
\begin{align}\label{g43}
\nonumber {A_1} &= 2{k_1}{\mu _1}\lambda {\rho _2}(1 - {\rho _0})(1 - \rho _0^N)(1 - \rho _0^{N + 1})(\rho _0^{N + 1} + \rho _0^N)(1 - {\rho _1}), \\
\nonumber {A_2} &= {k_1}{\mu _1}\rho _0^N{\rho _2}{(1 - {\rho _0})^2}(1 - \rho _0^{N + 1}),\\
\nonumber {A_3} &= {k_2}{\mu _1}{\rho _2}(1 - {\rho _0})\rho _0^N(1 - \rho _0^N),\\
\nonumber {A_4} &= 2{k_2}{\rho _1}{\rho _2}{(1 - {\rho _0})^2}{\rho _0}^{2N - 1}(1 - \rho _0^{N + 1}),\\
{B_1} &= {\rho _1}(\rho _0^{N + 1} + \rho _0^N) + (1 - \rho _0^{N + 1})(1 - {\rho _1}).
\end{align}
So we can obtain ${S'_M}(q) <0$ for $q\in [0, 1]$.

To sum up the above discussion, $S'(q) = {S_1}^\prime (q) + {S_2}^\prime (q) + {S_M}^\prime (q)<0$ for $q \in [0, 1]$, so $S(q)$ is a strictly decreasing function for $q \in [0, 1]$. Hence, we obtain the maximum of $S(q)$ in ${q^*} = \min (\frac{{dS(q)}}{{dq}}\left| {_{q = {q^*}}} \right. = 0, 1)$, so the number and the type of the socially optimal strategy is ${q^*} = \min (\frac{{dS(q)}}{{dq}}\left| {_{q = {q^*}}} \right. = 0, 1)$. \QED
\section{Observable case}
\label{sec:4}
For observable case, an arriving passenger can obtain the information of the queues of both passenger and the taxis. As before, if arriving passengers find that there are only taxis queueing (i.e. $N(t)<0$), they will join the system directly. If the arriving passengers find that the passenger queues length is $n$, i.e. $N(t)=n$. A threshold strategy will be employed, when $N(t)< n_{e}$, the arriving passengers join the system; when $N(t)= n_{e}$, the arriving passengers will balk.
\subsection{Equilibrium passenger threshold strategy}
For observable case, the service time of an individual passenger is $1/{\mu _2}$, and the matching time of an individual passenger is $k_{2}$. The utility function of an arriving passenger observes that the system state is $N(t)=n$ (i.e. the length of the passengers is $n$) is represented by $U(n)$, $U(n)$ has the following expression
\begin{equation}\label{g44}
U(n) = R - P - {C_P}\frac{{n + 1}}{{{\mu _2}}} - {C_{M,P}}{k_2}.
\end{equation}
In fact, for the service completion of an arriving passenger, he must wait for $n+1$ (including his own service time) service times and his own matching time $k_{2}$.

As mentioned above, if $N(t)< n_{e}$, the choice of arriving passengers is to join the system; if $N(t)= n_{e}$, the best response to an arriving passenger is balking. The value of threshold $n_{e}$ is given by the following theorem \ref{t535}.
\begin{theorem} \label{t535}
For observable case, equilibrium threshold strategy $n_{e}$ of a arriving passenger, upon joining the passenger-taxi double-ended queue with impatient passenger and two-point matching time is given
\begin{equation}\label{g45}
{n_e} = \left\lfloor {\frac{{(R - P - {C_{M,P}}{k_2}){\mu _2}}}{{{C_P}}}} \right\rfloor ,
\end{equation}
where $\left\lfloor x \right\rfloor $ is the largest integer not exceeding $x$.
\end{theorem}
{\bf Proof}~~
By the utility function of passenger of (\ref{g44}), the threshold $n_{e}$ needs to satisfy the following inequalities
\begin{equation}\label{g46}
\begin{cases}
U(n_{e}-1) \geq 0;\\
U(n_{e}) <0.
\end{cases}
\end{equation}
$n_{e}$ is a positive integer and by (\ref{g46}), we can get ${n_e} = \left\lfloor {\frac{{(R - P - {k_2}){\mu _2}}}{{{C_P}}}} \right\rfloor$. \QED
\subsection{Social optimization}
\begin{figure} [ht]
\centering
\includegraphics[width=12cm]{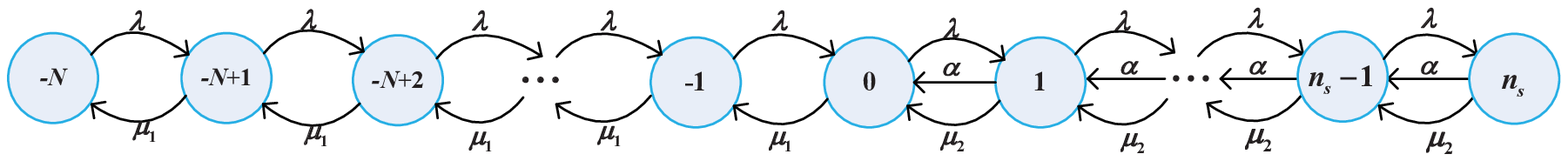}
\caption{The diagram of state transition for the observable case of the passenger-taxi double-ended queue with zero matching time.}
\label{Fig:g4}
\end{figure}

For observable case, the social welfare optimal strategy is threshold-type, which means that there exists a unique $n^{*}$ such that maximizes social welfare. So the state space is $\Omega _{o}= \{ -N, -N+1, -N+2, \cdots, -1, 0, 1, 2, \cdots, n_{s} \}$. The state transition diagram for the observable case with zero matching time is shown in Fig. \ref{Fig:g4}. According to the definition of the partially observable case, let ${\rho _0}{\rm{ = }}\frac{\lambda }{{{\mu _1}}}$, ${\rho _1}{\rm{ = }}\frac{{\lambda q}}{{\alpha  + {\mu _2}}}$ and ${\rho _2}{\rm{ = }}\frac{\lambda }{{\alpha  + {\mu _2}}}$. Let $\pi _n^o$ be the stationary probability of state $n \in {\Omega _o}$ for the observable case with zero matching time. The specific stationary probability is given by the following theorem \ref{t536}.
\begin{theorem} \label{t536}
For the observable case, the stationary probability of the passenger-taxi double-ended
queue with impatient passengers and zero matching time is given by
\begin{equation}\label{g49}
\pi _{ n}^o =
\begin{cases}
\pi _{ - N}^o{\rho _0}^{n + N}, &- N \le n \le 0;\\
\pi _{ - N}^o{\rho _0}^N{\rho _2}^n, &1 \le n \le {n_s};
\end{cases}
\end{equation}
where
\begin{equation}\label{g50}
\pi _{ - N}^o = \frac{{(1 - {\rho _0})(1 - {\rho _2})}}{{1 - {\rho _2} - \rho _0^{N + 1} + \rho _0^N{\rho _2} - \rho _0^N\rho _2^{{n_s} + 1} + \rho _0^{N + 1}\rho _2^{{n_s} + 1}}}.
\end{equation}
\end{theorem}
{\bf Proof}~~
According to Fig. \ref{Fig:g4}, the following equilibrium equations are obtained
\begin{align}\label{g51}
\lambda \pi _{ - N}^o &= {\mu _1}\pi _{ - N + 1}^o,\\
\label{g52}
(\lambda  + {\mu _1})\pi _n^o &= \lambda \pi _{n - 1}^o + {\mu _1}_{n + 1}^o, ~~- N + 1 \le n \le 1,\\
\label{g53}
(\lambda  + {\mu _1})\pi _0^o &= \lambda \pi _{ - 1}^o + (\alpha  + {\mu _2})\pi _1^o,\\
\label{g55}
(\lambda  + \alpha  + {\mu _2})\pi _n^o{\rm{ }} &= \lambda \pi _{n - 1}^o + (\alpha  + {\mu _2})\pi _{n + 1}^o,~~1 \le n \le {n_s}.
\end{align}
From (\ref{g51})-(\ref{g55}) and after some algebraic calculations, we can get (\ref{g49}). $\pi _{ - N}^o $ can be obtained by the normalization condition $\sum\limits_{n =  - N}^{{n_s}} {{\pi _n} = 1}  $. \QED

By (\ref{g49}), for the observable case, the effective arrival rate of passengers is
\begin{equation}\label{g56}
\lambda _P^e{\rm{ }} = \lambda \sum\limits_{n =  - N}^{ - 1} {{\pi _n} + \lambda \sum\limits_{n = 0}^{{n_s} - 1} {{\pi _n}} }  = \lambda \pi _{ - N}^o \left [\frac{{1 - \rho _0^N}}{{1 - {\rho _0}}} + \frac{{\rho _0^N(1 - \rho _2^{{n_s}})}}{{1 - {\rho _2}}} \right ],
\end{equation}
and the effective arrival rate of taxis is
\begin{equation}\label{g57}
\lambda _T^e{\rm{ }} = {\mu _1}\sum\limits_{n =  - N + 1}^0 {{\pi _n} + {\mu _2}\sum\limits_{n = 1}^{{n_s}} {{\pi _n}} }  = \pi _{ - N}^o \left [\frac{{\lambda (1 - \rho _0^N)}}{{1 - {\rho _0}}} + \frac{{{\mu _2}{\rho _2}\rho _0^N(1 - \rho _2^{{n_s}})}}{{1 - {\rho _1}}} \right ].
\end{equation}

For the observable case, denoting $E(L_P^o)$ and $E(L_T^o)$ as expected queue lengths for passenger and taxi, respectively, which have the following forms
\begin{equation}\label{g58}
E(L_P^o) = \sum\limits_{n = 0}^{{n_s}} {n{\pi _n} = \pi _{ - N}^o\rho _0^N\sum\limits_{n = 0}^{{n_s}} {n\rho _2^n} }  = \pi _{ - N}^o \rho _0^N\frac{{{\rho _2} + \rho _2^{{n_s} + 1}({n_s}{\rho _2} - {n_s} - 1)}}{{{{(1 - {\rho _2})}^2}}},
\end{equation}
\begin{equation}\label{g59}
E(L_T^o) = \sum\limits_{n =  - N}^0 { - n{\pi _n} = \pi _{ - N}^o\rho _0^N\sum\limits_{n = 0}^N {n{{\left(\frac{1}{{{\rho _0}}}\right)}^n}} }  = \frac{{\pi _{ - N}^o (\rho _0^{N + 1} - {\rho _0} + N - N{\rho _0})}}{{{{(1 - {\rho _0})}^2}}}.
\end{equation}

Similar to the discussion of lemma \ref{lem1}, the expected matching time for two-mass point distribution of the observable case is
\begin{equation}\label{g60}
E({M_o}) = \sum\limits_{n =  - N}^{{n_s}} {E(M\left| {N(t)} \right.)\pi _n^o = } \pi _{ - N}^o \left [{k_1} + {k_1}\frac{{{\rho _0}(1 - \rho _0^N)}}{{1 - {\rho _0}}} + {k_2}\rho _0^N\frac{{{\rho _2}(1 - \rho _2^{{n_s}})}}{{1 - {\rho _2}}} \right ].
\end{equation}

By the same method of (\ref{g25}), the social welfare function of the observable case is as follows
\begin{equation}\label{g61}
S(n) = \lambda _P^e(R - P - {C_P}E(W_P^o) - {C_{M,P}}E({M_o})) + \lambda _T^e(P - {C_T}E(W_T^o) - {C_{M,T}}E({M_o})),
\end{equation}
by using Little's law, $E(W_P^o) = E(L_P^o)/\lambda _P^e$ and $E(W_T^o) = E(L_T^o)/\lambda _T^e$. So (\ref{g61}) is equivalent to
\begin{equation}\label{g62}
S(n) = \lambda _P^e(R - P) + \lambda _T^eP - {C_P}E(L_P^o) - {C_T}E(L_T^o) - E({M_o})({C_{M,P}}\lambda _P^e + {C_{M,T}}\lambda _T^e).
\end{equation}

The optimal threshold of social welfare makes social welfare reach the maximum at some point, which is defined by $n^{*}$, and the optimal threshold $n^{*}$ needs to satisfy the following inequalities
\begin{equation}\label{g63}
\begin{cases}
S({n^*}) - S({n^*} + 1) \ge 0;\\
S({n^*}) - S({n^*} - 1) \ge 0.
\end{cases}
\end{equation}
By calculation, (\ref{g63}) is equivalent to
\begin{align}\label{g64}
\nonumber &{C_P}{\rho _2}[{n^*}(1 - {\rho _2}){A_5} - \rho _0^N{\rho _2}(1 - {\rho _0})(1 - \rho _2^{{n^*}})] + {E_2}\rho _2^{2{n^*}} \le M_{1}\\
&\le {C_P}{\rho _2}[({n^*} + 1)(1 - {\rho _2}){A_5} - \rho _0^N{\rho _2}(1 - {\rho _0})(1 - \rho _2^{{n^*} + 1})] + {E_2}\rho _2^{2({n^*} + 1)},
\end{align}
where
\begin{align}\label{g65}
\nonumber M_{1} &= (R - P)\lambda {A_5}{(1 - {\rho _2})^2} + {E_3}{\rho _2}({\rho _2} - 1) - {E_1} \Big[(1 - {\rho _2}){\rho _2}\\
&+ \frac{{{A_5}(1 - {\rho _2})}}{{\rho _0^N({\rho _0} - 1)}}\Big]+ {E_2}{\rho _2} + \frac{{2{A_5}{D_2}}}{{\rho _0^N({\rho _0} - 1)}} - \frac{{{A_5}{D_2}(1 + {\rho _2})}}{{\rho _0^N({\rho _0} - 1)}},
\end{align}
\begin{equation*}
{A_5} = 1 - {\rho _2} - \rho _0^{N + 1} + \rho _0^N{\rho _2},
\end{equation*}
\begin{align*}
{E_1} &= {\mu _2}\rho _0^N(1 - {\rho _0}){\rho _2}P - {C_{M,P}}[{k_2}\lambda \rho _0^N(1 - \rho _0^N){\rho _2} + {k_1}\lambda \rho _0^N(1 - {\rho _0}) + {k_1}\lambda \rho _0^{N + 1}(1 - \rho _0^N)]\\
&- {C_{M,T}}[{k_2}\lambda \rho _0^N(1 - \rho _0^N){\rho _2} + {k_1}{\mu _2}\rho _0^N{\rho _2}(1 - {\rho _0}) + {k_1}{\mu _2}\rho _0^{N + 1}{\rho _2}(1 - \rho _0^N)],
\end{align*}
\begin{equation*}
{E_2} = {C_{M,P}}{k_2}\lambda \rho _0^{2N}{\rho _2}(1 - {\rho _0}) + {C_{M,T}}{k_2}{\mu _2}\rho _0^{2N}\rho _2^2(1 - {\rho _0}),
\end{equation*}
\begin{align*}
{E_3} &= \lambda (1 - \rho _0^N)(1 - {\rho _2})P - {C_T}(1 - {\rho _2}) \left(N - \frac{{{\rho _0}(1 - \rho _0^N)}}{{(1 - {\rho _0})}}\right) - {C_{M,P}} \Big[ {k_1}\lambda (1 - \rho _0^N)(1 - {\rho _2}) \\
&+ \frac{{{k_1}\lambda {\rho _0}{{(1 - \rho _0^N)}^2}(1 - {\rho _2})}}{{(1 - {\rho _0})}}\Big] - {C_{M,T}}\Big[{k_1}\lambda (1 - \rho _0^N)(1 - {\rho _2}) + \frac{{{k_1}\lambda {\rho _0}{{(1 - \rho _0^N)}^2}(1 - {\rho _2})}}{{(1 - {\rho _0})}}\Big].
\end{align*}
Let
\begin{equation}\label{g66}
g(x) = {C_P}{\rho _2}{A_5}(1 - {\rho _2})x - {C_P}{\rho _2}\rho _0^N{\rho _2}(1 - {\rho _0})(1 - \rho _2^x) + {E_2}\rho _2^{2x},
\end{equation}
since for threshold strategy, which needs $n \geq 1$, so $x \geq 1$. The following theorem \ref{t537} gives the monotonicity of $g(x)$.
\begin{theorem} \label{t537}
The monotonicity of function $g(x)$ is related to the values of $\rho_{0}$ and $\rho_{2}$, specifically, it has the following cases
\begin{description}
\item[I:] If $0< \rho_{2} < \rho_{0} <1$, $g(x)$ is a decreasing function in $x \in [1, \infty ]$.
\item[II:] If $\rho_{0}>1$ and $\rho_{2} \in (0, 1)$, $g(x)$ is a increasing function in $x \in [1, \infty ]$.
\item[III:] If $\rho_{0} > \rho_{2} >1$, $g(x)$ is a decreasing function in $x \in [1, \infty ]$.
\end{description}
\end{theorem}
{\bf Proof}~~
The first order derivative of $g(x)$ with respect to $x$ is
\begin{equation}\label{g67}
f({\rho _2}) = g'(x) = {C_P}{\rho _2}{A_5}(1 - {\rho _2}) + {C_P}\rho _0^N(1 - {\rho _0})\ln {\rho _2}\rho _2^{x + 2} + 2{E_2}\ln {\rho _2}\rho _2^{2x},
\end{equation}
where $A_{5}$ and $E_{2}$ are given in (\ref{g65}). The first order derivative of $f(\rho_{2})$ with respect to $\rho_{2}$ is
\begin{align}\label{g68}
\nonumber f'({\rho _2}) &= {C_P}((1 - \rho _0^{N + 1}) + 3\rho _2^2(1 - \rho _0^N) + 2{\rho _2}(1 - \rho _0^N) + 2{\rho _2}(1 - \rho _0^{N + 1})) \\
\nonumber &+ 2(1 - {\rho _0}){k_2}\rho _0^{2N}\rho _2^{2x}({C_{M,P}}\lambda (1 + \ln {\rho _2} + x\ln \rho _2^2) + {C_{M,T}}{\mu _2}(1 + \ln \rho _2^2 + x\ln \rho _2^2))\\
&+ {C_P}(1 - {\rho _0})\rho _0^N\rho _2^{x + 1}(1 + (x + 2)\ln {\rho _2}).
\end{align}
Obviously, ${\rho _0}{\rm{ = }}\frac{\lambda }{{{\mu _1}}} > {\rho _2}{\rm{ = }}\frac{\lambda }{{\alpha  + {\mu _2}}}$, according to the values of $\rho_{0}$ and $\rho_{2}$, the monotonicity of the function $g(x)$ is as follows\\
(1) If $\rho_{0}<1$, from (\ref{g68}), which produces $f'({\rho _2})>0$, so $f({\rho _2})$ is an increasing function with respect to $\rho_{2}$. Then
\begin{equation}\label{g69}
f({\rho _2}) = g'(x) < f(1) = 0,~~{\rho _2} \in (0,1),
\end{equation}
i.e., when $0< \rho_{2} < \rho_{0} <1$, $g(x)$ is a decreasing function in $x \in [1, \infty ]$.\\
(2) If $\rho_{0}>1$, from (\ref{g68}), which produces $f'({\rho _2})<0$, so $f({\rho _2})$ is a decreasing function with respect to $\rho_{2}$. Then
\begin{equation}\label{g70}
f({\rho _2}) = g'(x) > f(1) = 0,~~{\rho _2} \in (0,1),
\end{equation}
i.e., when $\rho_{0}>1$ and $\rho_{2} \in (0, 1)$, $g(x)$ is a increasing function in $x \in [1, \infty ]$. on the other hand,
\begin{equation}\label{g71}
f({\rho _2}) = g'(x) < f(1) = 0,~~{\rho _2} \in (1, + \infty ),
\end{equation}
i.e., $\rho_{0} > \rho_{2} >1$, $g(x)$ is a decreasing function in $x \in [1, \infty ]$. \QED

Let
\begin{equation}\label{g73}
{C_P}{\rho _2}[{n^*}(1 - {\rho _2}){A_5} - \rho _0^N{\rho _2}(1 - {\rho _0})(1 - \rho _2^{{n^*}})] + {E_2}\rho _2^{2{n^*}}{\rm{ = }}M_{1},
\end{equation}
where $M_{1}$ is given in (\ref{g65}). Obviously, the solution of the equation (\ref{g73})
is the socially optimal threshold $n^{*}$. The following theorem \ref{t538} clearly shows the existence of the root of the equation (\ref{g73}), and the specific proof is no longer described.
\begin{theorem} \label{t538}
The monotonicity of function $g(x)$ is related to the values of $\rho_{0}$ and $\rho_{2}$, specially, it has the following cases
\begin{description}
\item[I:] If $\rho_{0}=1$ or $\rho_{2}=1$, the equation of (\ref{g73}) has no solution.
\item[II:] For $0< \rho_{2} < \rho_{0} <1$,
\begin{description}
\item[(1)] If $g(1)<M_{1}$, the equation of (\ref{g73}) has no solution.
\item[(2)] If $g(1)>M_{1}$, the equation of (\ref{g73}) has unique solution $n^{*} >1$.
\item[(3)] If $g(1)=M_{1}$, the equation of (\ref{g73}) has unique solution $n^{*} =1$.
\end{description}

\item[III:] If $\rho_{0}>1$ and $\rho_{2} \in (0, 1)$,
\begin{description}
\item[(1)] If $g(1)<M_{1}$, the equation of (\ref{g73}) has unique solution $n^{*} >1$.
\item[(2)] If $g(1)>M_{1}$, the equation of (\ref{g73}) has no solution.
\item[(3)] If $g(1)=M_{1}$, the equation of (\ref{g73}) has unique solution $n^{*} =1$.
\end{description}

\item[IIII:] If $\rho_{0} > \rho_{2} >1$,
\begin{description}
\item[(1)] If $g(1)<M_{1}$, the equation of (\ref{g73}) has no solution.
\item[(2)] If $g(1)>M_{1}$, the equation of (\ref{g73}) has unique solution $n^{*} >1$.
\item[(3)] If $g(1)=M_{1}$, the equation of (\ref{g73}) has unique solution $n^{*} =1$.
\end{description}
\end{description}
\end{theorem}
\section{Numerical examples}
\label{sec:5}
From the previous theoretical results, we know that the system of the passenger-taxi double-ended queue with impatient passengers and two-point matching time constructed in this paper contains a large number of parameters, the change of these parameters has a great influence on the equilibrium joining probability $q_{e}$ and socially optimal joining probability $q^{*}$ of the partially observable case, as well as the equilibrium threshold $n_{e}$ and socially optimal threshold $n^{*}$ of the observable case. In addition, the change trend of optimal social welfare on the same parameter is also presented under two information levels, which has some guiding significance for the manager to reveal or hide system information.

\begin{figure}[ht]
\centering
\subfigure[]{
\includegraphics[width=6.5cm]{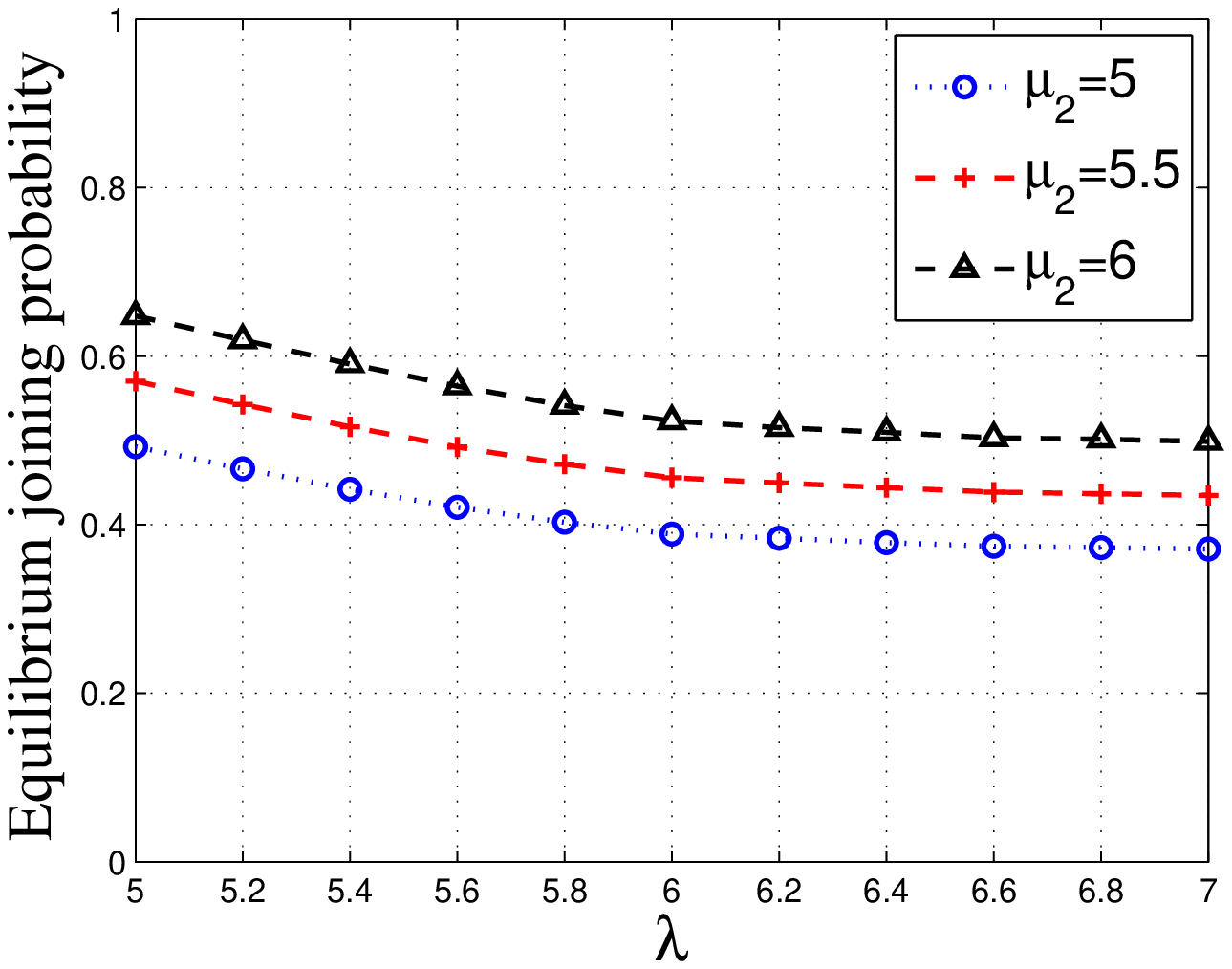}
}
\quad
\subfigure[]{
\includegraphics[width=6.5cm]{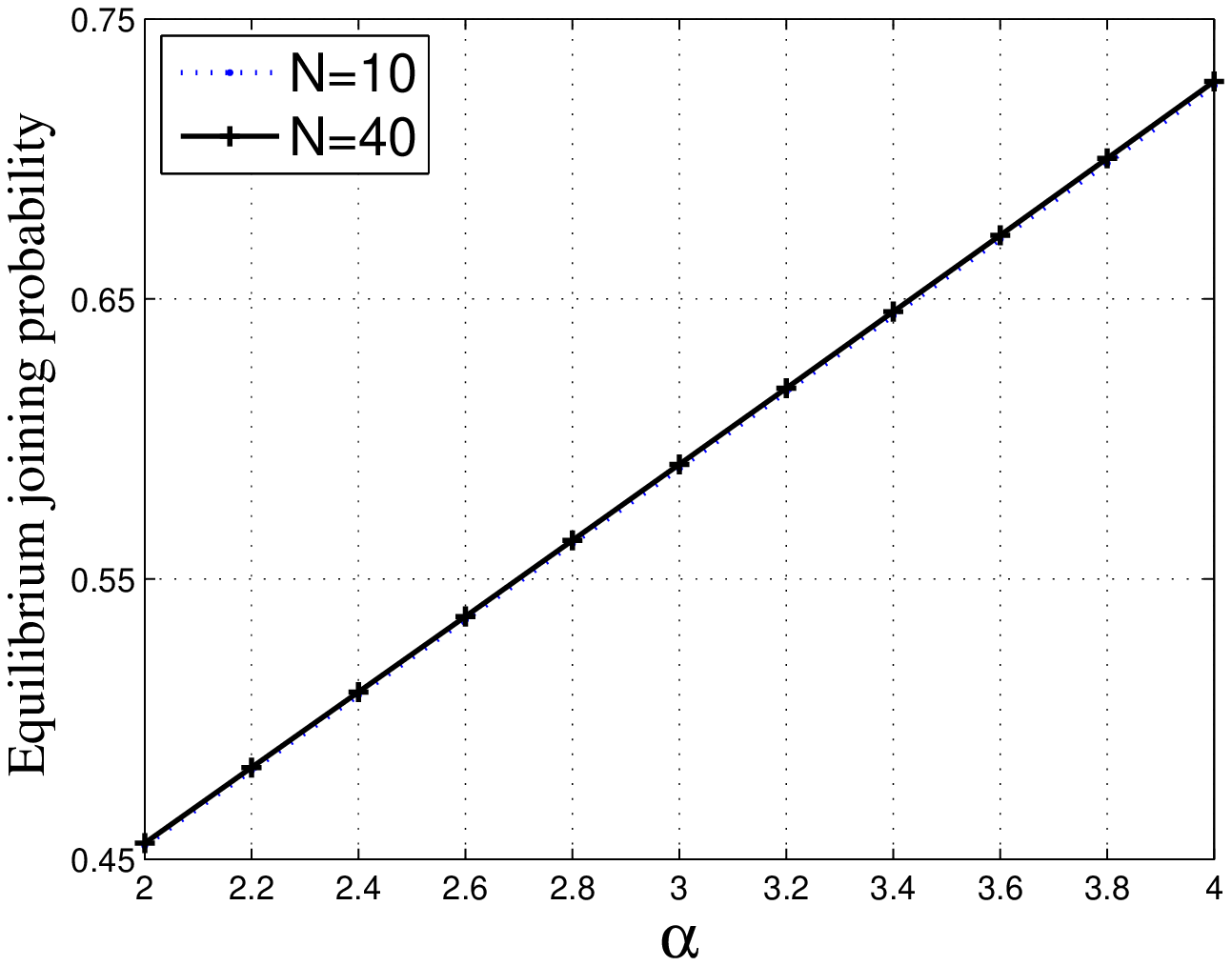}
}
\caption{(a) Equilibrium joining probability $q_{e}$ with respect to $\lambda$ and $\mu_{2}$ when $R=16$, $P=6$, $C_{P}=4$, $\alpha=2$, $C_{M,P}=3$, $\mu_{1}=4$, $N=20$, $k_{1}=3$, $k_{2}=5$.  (b) Equilibrium joining probability $q_{e}$ with respect to $\alpha$ and $N$ when $R=16$, $P=6$, $C_{P}=4$, $\lambda=6$, $\mu_{2}=5.5$,  $C_{M,P}=3$, $\mu_{1}=4$, $k_{1}=3$, $k_{2}=5$.}
\label{Fig:g5}
\end{figure}

Firstly, we explore the effects of the potential arrival rate $\lambda$ of passengers, taxi arrival rate $\mu_{2}$, the passenger's impatience rate $\alpha$ and taxi queue capacity $N$ on the equilibrium joining probability $q_{e}$ in the partially observable case. Extensive numerical experiments for system parameter selection are performed. It has been found that the numerical changes presented in this section are generally insensitive to the choice of parameters. Therefore, we only present typical numerical scenarios to illustrate our research results. These typical numerical results are shown in Fig. \ref{Fig:g5}. \\
1. For the arriving passengers of the partially observable case, the equilibrium joining probability $q_{e}$ decreases with the increase of the potential arrival rate $\lambda$ (see Fig. \ref{Fig:g5} (a)), which is because with the increase of $\lambda$, the system becomes more and more congested, which leads to the decrease of the joining probability $q_{e}$ of arriving passengers. On the contrary, the increase of $\mu_{2}$ makes the system more fluent, which also promotes the arriving passengers to be more willing to join the system, so $q_{e}$ is increases as $\mu_{2}$ increases (see Fig. \ref{Fig:g5} (a)).\\
2. From \ref{Fig:g5} (b), it seems like the $N=10$ and $N=40$ lines are almost the same. From the figure, it is not clear that increasing $N$ has an effect on $q_e$.

\begin{figure}[ht]
\centering
\subfigure[]{
\includegraphics[width=6.5cm]{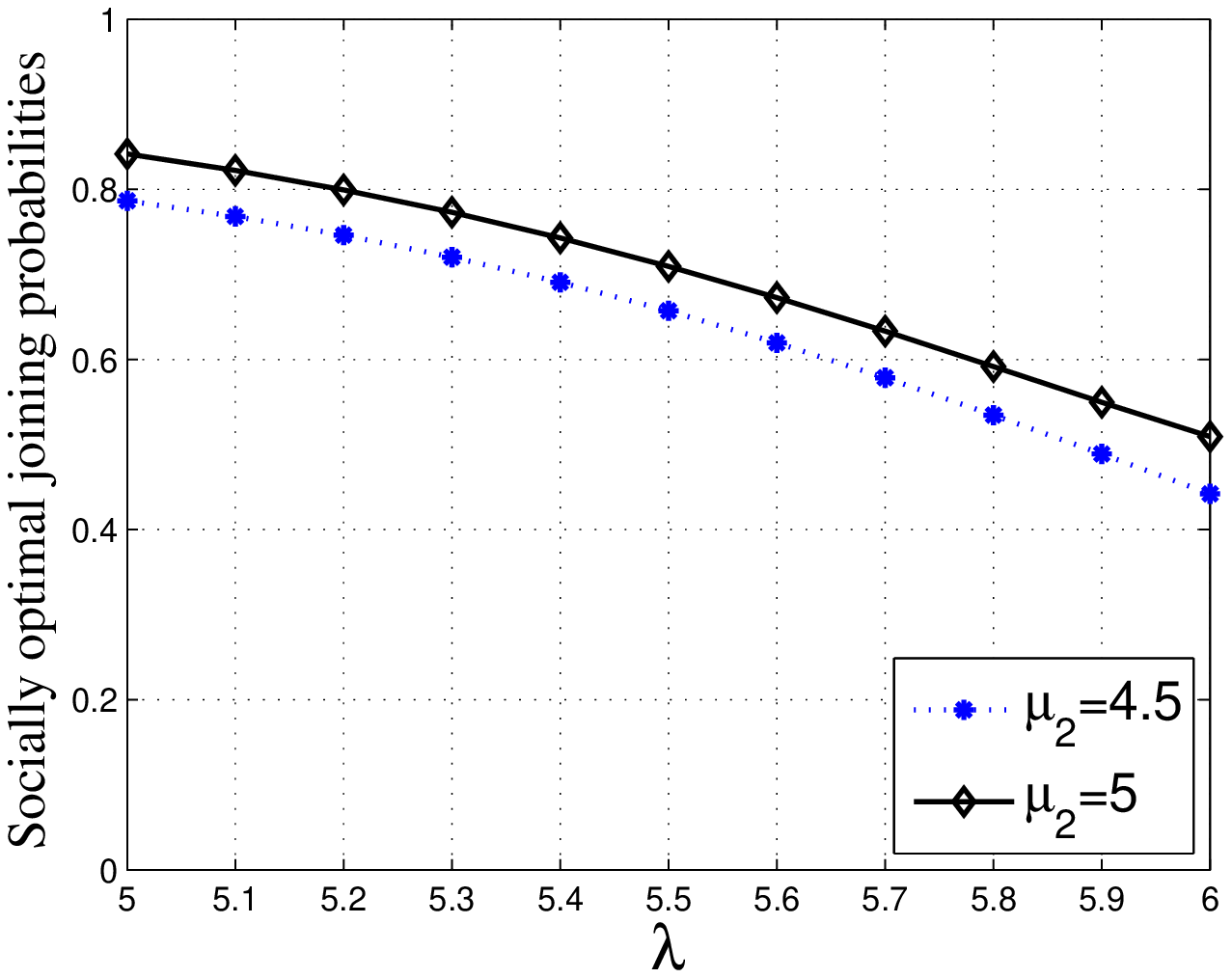}
}
\quad
\subfigure[]{
\includegraphics[width=6.5cm]{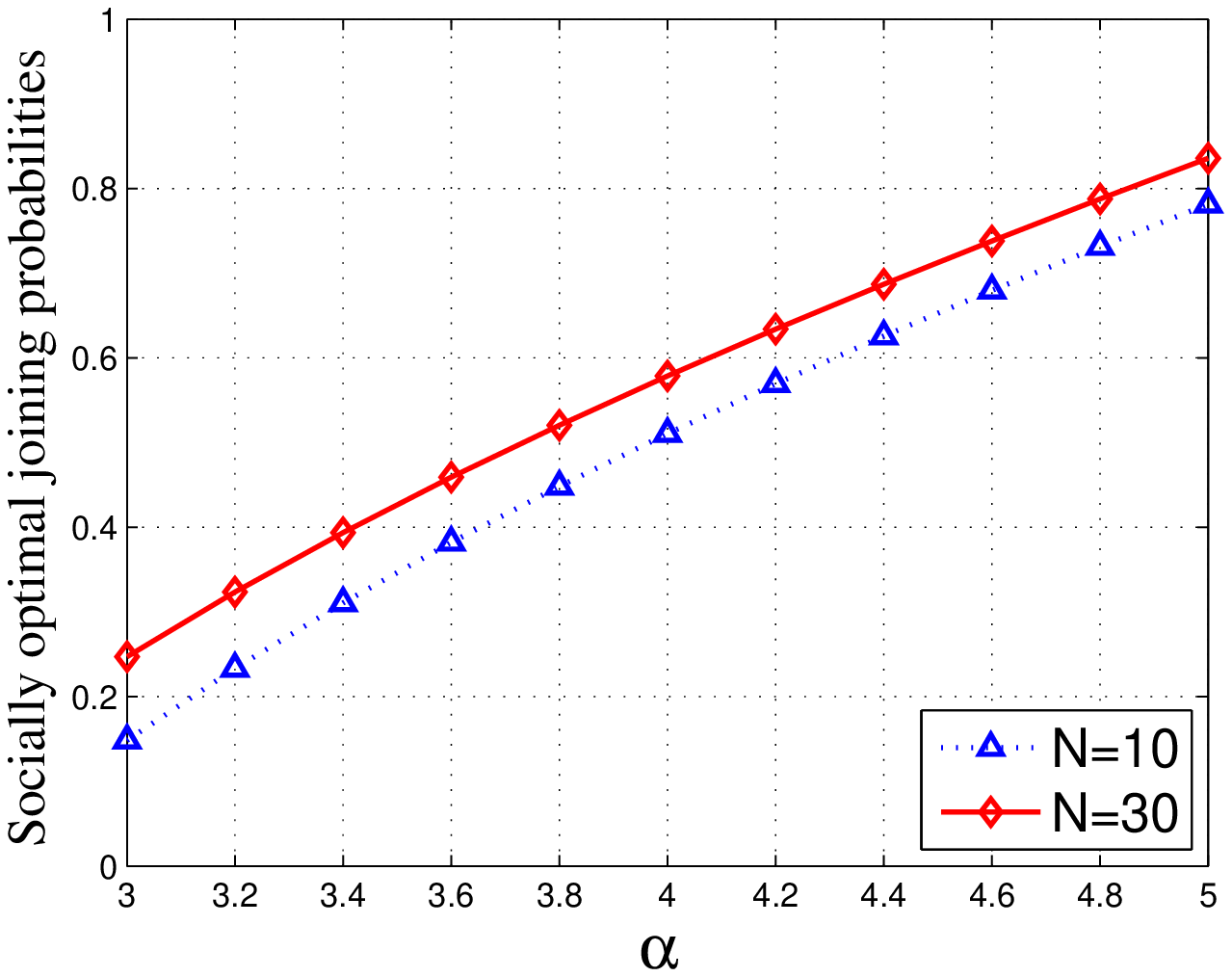}
}
\caption{(a) Socially optimal joining probability $q^*$ with respect to $\lambda$ and $\mu_{2}$ when $R=20$, $P=6$, $C_{P}=4$, $C_{T}=3$, $\alpha=4$, $C_{M,P}=3$, $C_{M,T}=3$, $\mu_{1}=4$, $N=30$, $k_{1}=3$, $k_{2}=5$.  (b) Socially optimal joining probability $q^*$ with respect to $\alpha$ and $N$ when $R=20$, $P=6$, $C_{P}=4$, $C_{T}=3$, $\mu_{2}=4.5$, $\lambda=5.3$,  $C_{M,P}=3$, $C_{M,T}=3$, $\mu_{1}=4$, $k_{1}=3$, $k_{2}=5$.}
\label{Fig:g6}
\end{figure}

Secondly, we consider the influence of the potential arrival rate $\lambda$ of passengers, taxi arrival rate $\mu_{2}$, the passenger's impatience rate $\alpha$ and maximum taxi arrival value $N$ on the socially optimal joining probability $q^{*}$ in the partially observable case. These typical numerical results are shown in Fig. \ref{Fig:g6}. \\
1. For the arriving passengers of the observable case, the socially optimal joining probability $q^*$ decreases with the increase of arrival rate $\lambda$ (see Fig. \ref{Fig:g6} (a)), because the arriving passengers can observe both the passengers and the taxi, the congestion of the system has a negative impact on the arriving passengers. At the same time, the increase in taxi arrival rate $\mu_{2}$ accelerates the operation of the system, which encourages arriving customers to join the system, so $q^*$ is increasing as $\mu_{2}$ increases (see Fig. \ref{Fig:g6} (a)).\\
2. From the perspective of system operation, the loss of impatient customers reduces the congestion of the system, which encourages the late-arriving passengers to join the system, so $q^*$ is increasing as $\alpha$ increases (see Fig. \ref{Fig:g6} (b)). Similarly, the increase of capacity $N$ of taxis also accelerates the operation of passengers, so $q^*$ is increasing as $N$ increases (see Fig. \ref{Fig:g6} (b)).

\begin{figure}[h!]
\centering
\subfigure[]{
\includegraphics[width=6.5cm]{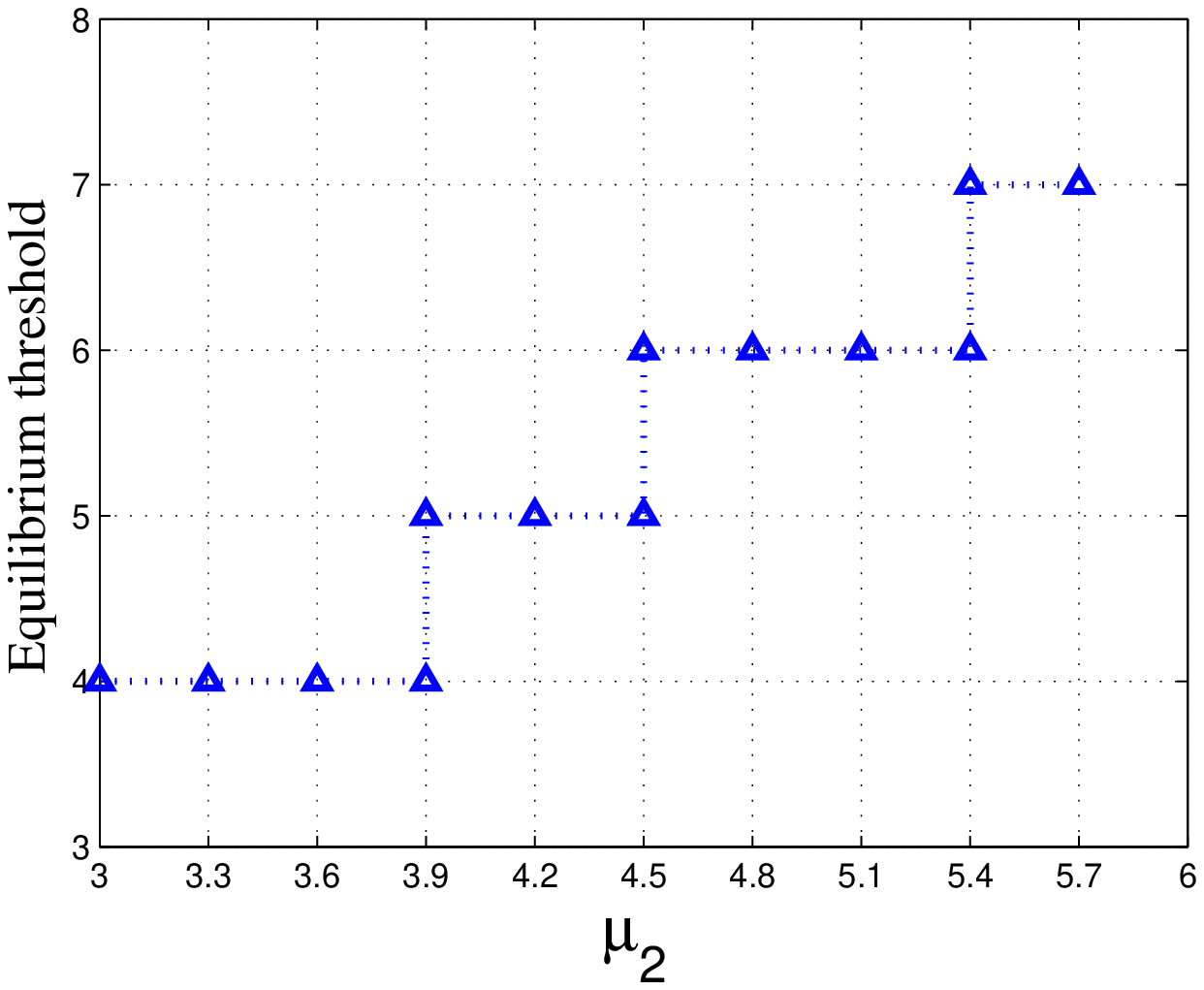}
}
\quad
\subfigure[]{
\includegraphics[width=6.5cm]{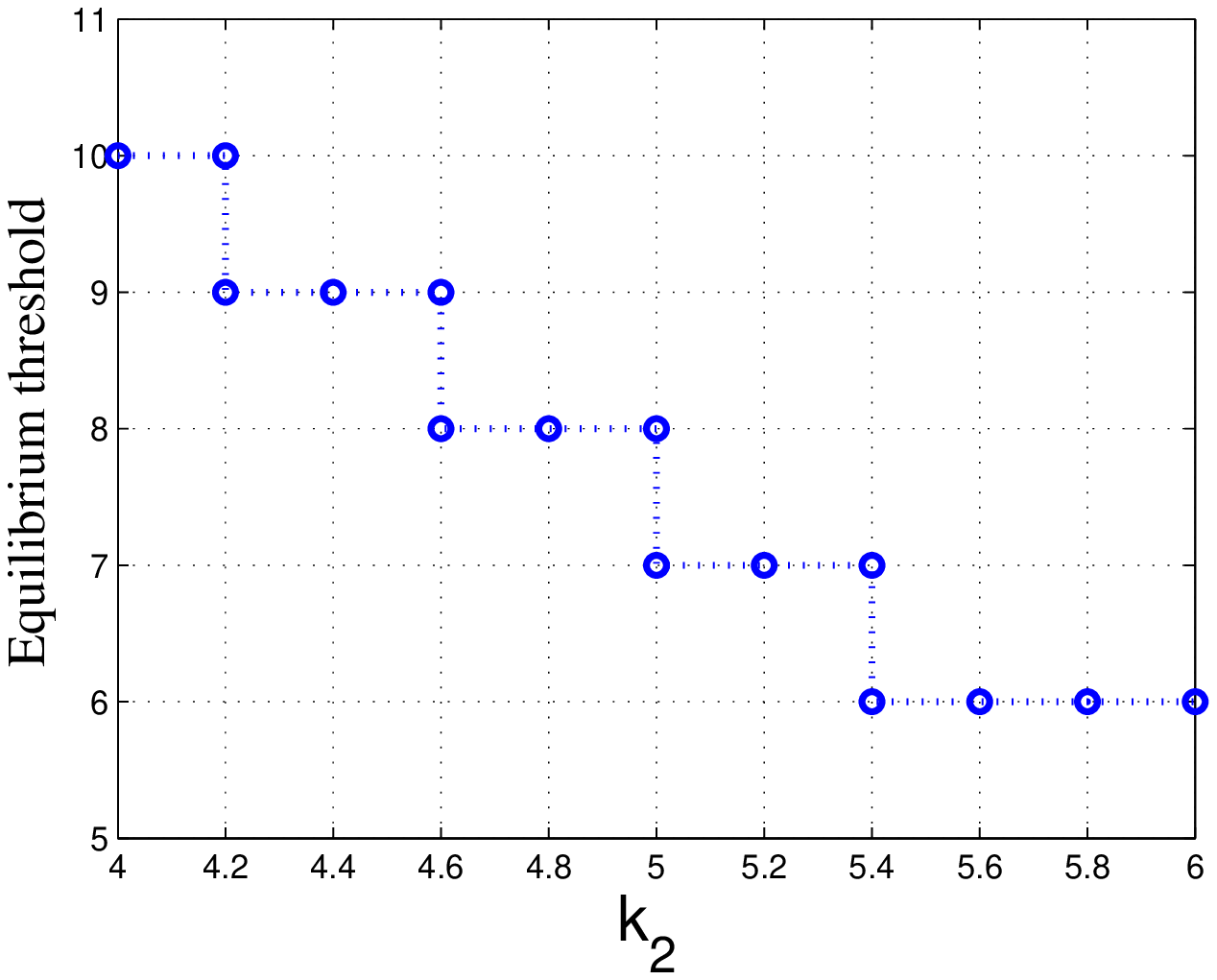}
}
\caption{(a) Equilibrium threshold $n_{e}$ with respect to $\mu_{2}$ when $R=15$, $P=6$, $C_{P}=3$, $k_{2}=5$.  (b) Equilibrium threshold $n_{e}$ with respect to $k_{2}$ when $R=15$, $P=6$, $C_{P}=3$, $\mu_{2}=6$.}
\label{Fig:g7}
\end{figure}

\begin{figure}[h!]
\centering
\subfigure[]{
\includegraphics[width=6.5cm]{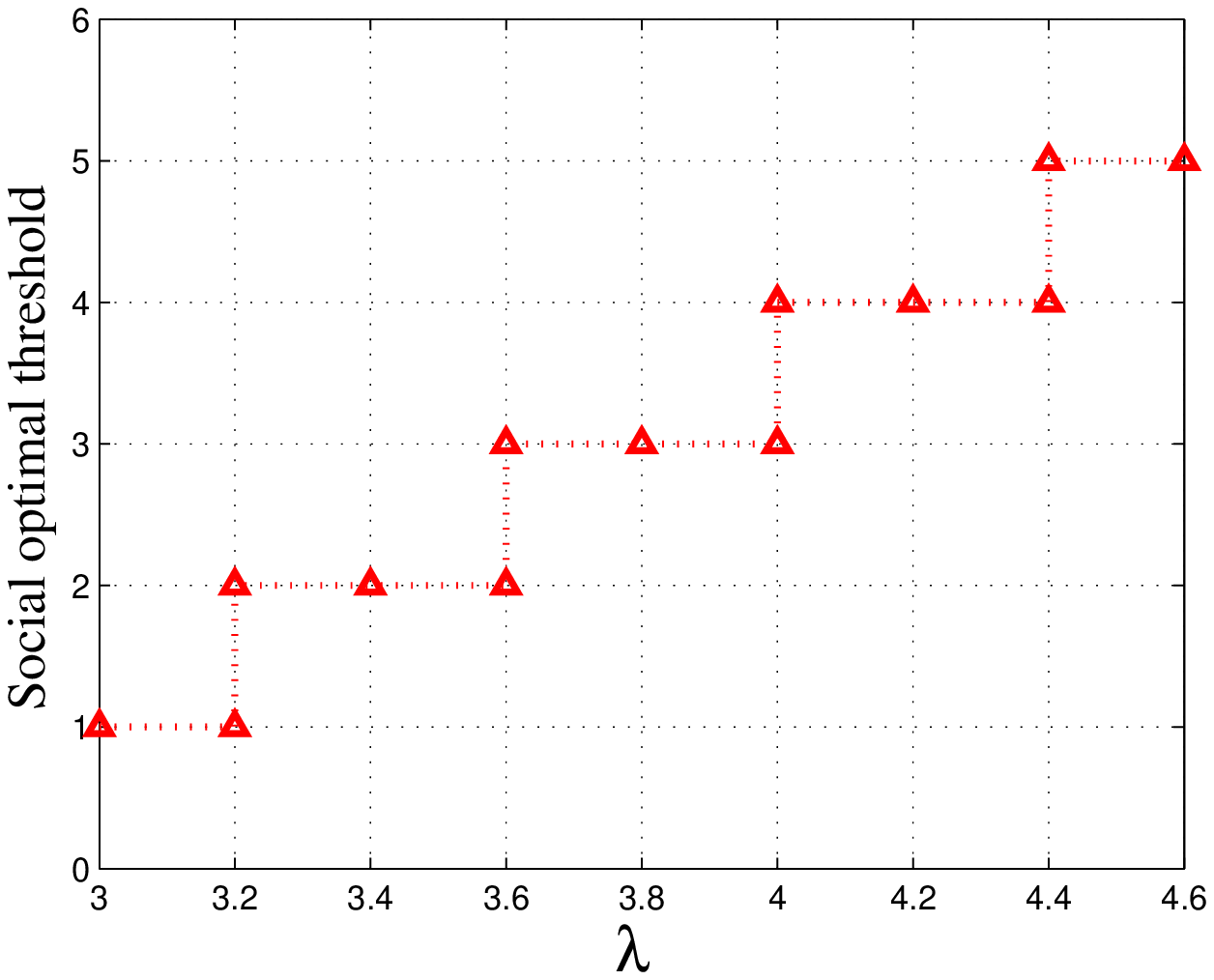}
}
\quad
\subfigure[]{
\includegraphics[width=6.5cm]{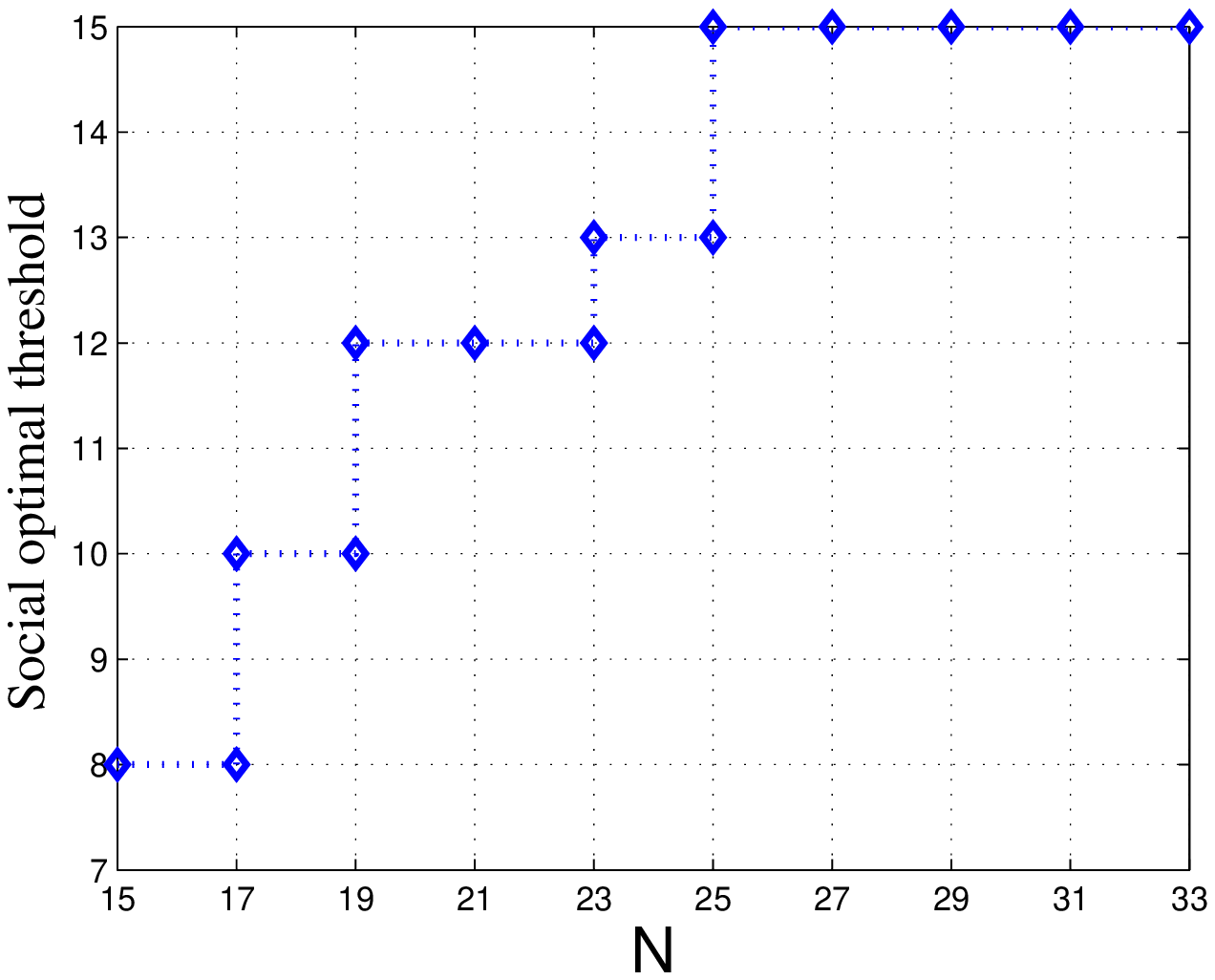}
}
\caption{(a) Social optimal threshold $n^{*}$ with respect to $\lambda$ when $R=20$, $P=6$, $C_{P}=4$, $C_{T}=3$, $\alpha=4$, $\mu_{2}=4.5$, $C_{M,P}=3$, $C_{M,T}=3$, $\mu_{1}=4$, $N=30$, $k_{1}=3$, $k_{2}=5$.  (b) Social optimal threshold $n^{*}$ with respect to $N$ when $R=20$, $P=6$, $C_{P}=4$, $C_{T}=3$, $\alpha=4$, $\mu_{2}=4.5$, $\lambda=5.3$, $C_{M,P}=3$, $C_{M,T}=3$, $\mu_{1}=4$, $N=30$, $k_{1}=3$, $k_{2}=5$.}
\label{Fig:g8}
\end{figure}

Thirdly, we study the influence of parameters on the equilibrium threshold $n_{e}$ and the socially optimal threshold $n^{*}$, which are shown in figures \ref{Fig:g7} and \ref{Fig:g8},  respectively.

$n_{e}$ is increasing as $\mu_{2}$ increases (see Fig. \ref{Fig:g7} (a)), which is because the increase of taxi arrival rate $\mu_{2}$ accelerates the transfer of passengers, and promotes the increase of social welfare and equilibrium threshold $n_{e}$. $n_{e}$ is decreasing as $k_{2}$ increases (see Fig. \ref{Fig:g7} (b)), the reason is that the increase of $k_{2}$ increases the cost of matching time of passengers, which will lead to a decrease in passengers joining the system.

For the socially optimal threshold $n^{*}$ of the observable case, the increase of $\lambda$ and $N$ both can promote its growth (see Fig. \ref{Fig:g8}), because these measures promote the growth of social welfare and equilibrium threshold $n_{e}$.

Finally, under the same parameters, we compare the optimal social welfare under two information levels (partially observable case and observable case). The numerical results show that the levels of information disclosure have a great effect on the increase in social welfare. The fully observable information level does not necessarily promote the growth of social benefit, and the concealment of the information level does not necessarily lead to a negative effect.
\begin{figure} [ht]
\centering
\includegraphics[width=6.5cm]{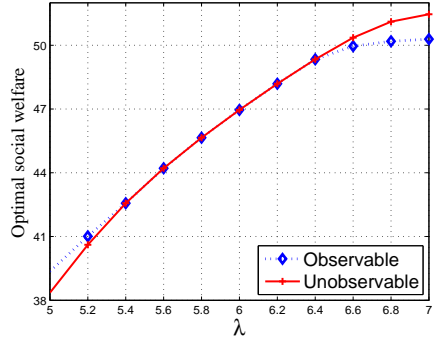}
\caption{Optimal social welfare under two different information levels with respect to $\lambda$ when $R=20$, $P=6$, $C_{P}=4$, $C_{T}=3$, $\alpha=4$, $\mu_{2}=4.5$, $C_{M,P}=3$, $C_{M,T}=3$, $\mu_{1}=4$, $N=30$, $k_{1}=3$, $k_{2}=5$.}
\label{Fig:g9}
\end{figure}
\begin{figure} [ht]
\centering
\includegraphics[width=6.5cm]{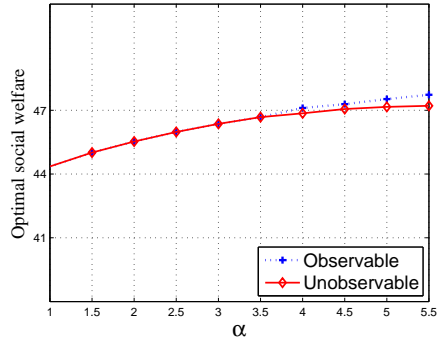}
\caption{Optimal social welfare under two different information levels with respect to $\alpha$ when $R=20$, $P=6$, $C_{P}=4$, $C_{T}=3$, $\mu_{2}=4.5$, $\lambda=5.3$, $C_{M,P}=3$, $C_{M,T}=3$, $\mu_{1}=4$, $N=30$, $k_{1}=3$, $k_{2}=5$, $N=20$.}
\label{Fig:g10}
\end{figure}

Fig. \ref{Fig:g9} shows the impact of passenger arrival rate $\lambda$ on optimal social welfare under two information levels, it contains a lot of numerical information, which has important guiding significance for when to reveal or hide the information. The specific explanation is as follows,\\
1. For a smaller arrival rate $\lambda$, the system is not crowded, so the arriving passengers are more willing to join the system. The arriving passengers in the observable case make effective use of this information, so the optimal social welfare of the observable case is greater than the optimal social welfare of the partially observable case.\\
2. For the moderate passenger arrival rate $\lambda$, the optimal social welfare under the two information levels is consistent.\\
3. For a large arrival rate $\lambda$, it causes the system to become crowded, and the arriving passengers of the observable case are more reluctant to join the system. At this time, the information hiding promotes more arriving passengers to join the system and further promotes the growth of the optimal social welfare.

Fig. \ref{Fig:g10} shows the impact of the impatience rate $\alpha$ of passengers on optimal social welfare under two information levels, the growth of impatience rate $\alpha$ can promote the optimal social welfare under the two information levels. For small and moderate impatience rate $\alpha$, the optimal social welfare under the two information levels is consistent. For the larger impatience rate $\alpha$, revealing this information can promote the arriving passengers to join the system, and further increase the optimal social welfare. Therefore, for a large impatience rate, the optimal social welfare of observable cases is greater than that of partially observable case.

\section{Conclusions }
\label{sec:6}

In this paper, we first derive the performance measures of the passenger-taxi double-ended queue with impatient passengers and zero matching time. Then, under two information levels, we use these performance measures to explore the equilibrium strategy and socially optimal strategy of the system with the two-point matching time. The theoretical results show that the passenger utility function of the partially observable case is monotonic. For the complex form of social welfare function of the partially observable case, we use split derivation. The equilibrium strategy and socially optimal strategy of the observable case are threshold-type. The model considered in this paper is more practical, and it is specifically manifested in several aspects. For different passenger queues, the taxis arrival rate is dynamically controlled, the capacity of the taxi is limited, and the waiting passengers have impatient behavior. In addition, different from the existing literature, the system considered in this paper fully considers the matching time between the passenger and the taxi. The consideration of this factor promotes the studied system's practically and complexity.

To visualizing the theoretical results, typical numerical experiments are presented. Some typical-numerical scenarios illustrate the influence of parameters on equilibrium strategy and socially optimal strategy under the two information levels. As well, the optimal social welfare of two information levels with the same parameters are compared, which shows that the disclosure or concealment of information can promote the growth of social welfare. This system with priority consideration is a future research topic that can reduce the waiting time of passengers.

\section*{Acknowledgements}
This work was supported in partial by The National Natural Science Foundation of China (No. 61773014), the Research Fund for the Postgraduate Research and Practice Innovation Program of Jiangsu Province (No. KYCX20\_0240), and the Natural Sciences and Engineering Research Council of Canada (NSERC).
\section*{Disclosure statement}
None of the authors have any competing interests in the manuscript.


\begin{thebibliography}{10}
\expandafter\ifx\csname url\endcsname\relax
  \def\url#1{\texttt{#1}}\fi
\expandafter\ifx\csname urlprefix\endcsname\relax\def\urlprefix{URL }\fi
\expandafter\ifx\csname href\endcsname\relax
  \def\href#1#2{#2} \def\path#1{#1}\fi

\bibitem{Curry1978A}
G.L. Curry, A.D. Vany, R.M. Feldman, A queueing model of airport passenger
  departures by taxi: Competition with a public transportation mode,
  Transportation Research 12~(2) (1978) 115--120.

\bibitem{2013Passos}
L.S. Passos, Z. Kokkinogenis, R. Rossetti, J. Gabriel, Multi-resolution simulation of taxi services on airport terminal\"s curbside, in: Itsc. 16th international ieee conference on intelligent transportation systems, the hague, Netherlands, 2013, pp. 2361--2366.

\bibitem{Conway2018Challenges}
A. Conway, C. Kamga, A. Yazici, A. Singhal, Challenges in managing centralized taxi dispatching at high-volume airports: Case study of John F. Kennedy International Airport, New York City, Transportation Research Record, 2300~(1) (2012) 83--90.

\bibitem{Yazici2016Modeling}
M. Yazici, C. Kamga, A. Singhal, Modeling taxi drivers' decisions for
  improving airport ground access: John f. kennedy airport case, Transportation
  Research Part A: Policy and Practice 91 (2016) 48--60.

\bibitem{Kendall1951Some}
D.G. Kendall, Some problems in the theory of queues, Journal of the Royal
  Statistical Society 13~(2) (1951) 151--185.

\bibitem{Dobbie1961Letter}
J.M. Dobbie, M.~James, A doubled-ended queuing problem of kendall, Operations
  Research 9~(5) (1961) 755--757.

\bibitem{wang2020equilibrium}
Z. Wang, L. Liu, Y. Shao, X. Chai, B. Chang, Equilibrium joining strategy in a
  batch transfer queuing system with gated policy, Methodology and Computing in
  Applied Probability 22~(1) (2020) 75--99.

\bibitem{giveen1963taxicab}
S.M. Giveen, A taxicab problem with time-dependent arrival rates, SIAM Review
  5~(2) (1963) 119--127.

\bibitem{Conolly2002Double}
W. Conolly, R. Parthasarathy, N. Selvaraju, Double-ended queues with
  impatience, Computers and Operations Research, 29~(14) (2002) 2053--2072.

\bibitem{Kashyap1965A}
B. Kashyap, A double ended queuing system with limited waiting space,
  Proceedings of the National Institute of Science of India 31 (1965) 559--570.

\bibitem{Naor1969The}
P. Naor, The Regulation of Queue Size by Levying Tolls, Econometrica 37~(1)
  (1969) 15--24.

\bibitem{Edelson1975Congestion}
N.M. Edelson, D.K. Hilderbrand, Congestion Tolls for Poisson Queuing
  Processes, Econometrica 43~(1) (1975) 81--92.

\bibitem{Jiang2019Tail}
T. Jiang, L. Liu, Q. Ye, X. Chai, Tail asymptotics for service systems with
  transfers of customers in an alternating environment, Operations Research
  Letters 47~(6) (2019) 473--477.

\bibitem{boudali2012optimal}
O. Boudali, A. Economou, Optimal and equilibrium balking strategies in the
  single server markovian queue with catastrophes, European Journal of
  Operational Research 218~(3) (2012) 708--715.

\bibitem{economou2013equilibrium}
A. Economou, A. Manou, Equilibrium balking strategies for a clearing queueing
  system in alternating environment, Annals of Operations Research 208~(1)
  (2013) 489--514.

\bibitem{burnetas2007equilibrium}
A. Burnetas, A. Economou, Equilibrium customer strategies in a single server
  markovian queue with setup times, Queueing Systems 56~(3-4) (2007) 213--228.


\bibitem{sun2017equilibrium}
W. Sun, S. Li, N. Tian, Equilibrium and optimal balking strategies of customers
  in unobservable queues with double adaptive working vacations, Quality
  Technology and Quantitative Management 14~(1) (2017) 94--113.

\bibitem{Bu2020Strategic}
Q. Bu, Y. Sun, X. Chai, L. Liu, Strategic behavior and social optimization in a clearing queueing system with N-policy and stochastic restarting scheme, Applied Mathematics and Computation 381 (2020) 125309.

\bibitem{Hassin2003To}
R. Hassin, M. Haviv,  To Queue or Not to Queue: Equilibrium Behavior in
  Queueing Systems, Kluwer Academic Publishers, Boston, 2003.

\bibitem{Shi2016Optimization}
Y. Shi, Z. Lian, Optimization and strategic behavior in a passenger-taxi
  service system, European Journal of Operational Research (2016) 1024--1032.

\bibitem{Wang2019Equilibrium}
Y. Wang, Z. Liu, Equilibrium and optimization in a double-ended queueing system
  with dynamic control, Journal of Advanced Transportation 2019~(2) (2019)
  1--13.

\bibitem{Kim2010Simulation}
W.K. Kim, K.P. Yoon, G. Mendoza, M. Sedaghat, Simulation model for extended
  double-ended queueing, Computers and Industrial Engineering 59~(2) (2010)
  209--219.

\bibitem{Shi2015Study}
Y. Shi, Z. Lian, W. Shang, Study of a passenger-taxi queueing system with
  nonzero matching time, in: International Conference on Service Systems and
  Service Management, Guangzhou, China, 2015, pp. 1--5.

\end{thebibliography}
%

\end{document}